\documentclass[reqno]{amsart}
\usepackage{amsmath}
\usepackage{amsthm}
\usepackage{amssymb}
\usepackage{epsfig}
\usepackage{psfrag}
\usepackage{graphicx}
\usepackage{pgf}
\newcommand{\disp}{\displaystyle}
\newcommand{\Rz}{\mathbb{R}}
\newcommand{\Tz}{\mathbb{T}}
\newcommand{\Zz}{\mathbb{Z}}

\newcommand{\Nz}{\mathbb{N}}

\newcommand{\dx}{{\rm d} x} 
\newcommand{\dr}{{\rm d} r} 
\newcommand{\ove}{\overline}
\newcommand{\haz}{\widehat}
\newcommand{\argmin}{{\rm arg\,min}}
\renewcommand{\sl}{|\partial \phi|}
\newcommand{\md}{|u'|}
\newcommand{\sto}{\stackrel{\sigma}{\to}}

\newcommand{\DD}{{\rm D}}

 \newtheorem{theorem}{Theorem}[section]

 \newtheorem{definition}[theorem]{Definition}
 \newtheorem{proposition}[theorem]{Proposition}

\begin{document}

\title[A new minimizing-movements scheme]{A new minimizing-movements
  scheme \\ for curves of maximal slope}
 
\author[U. Stefanelli]{Ulisse Stefanelli} 
\address[Ulisse Stefanelli]{Faculty of Mathematics, University of
  Vienna, Oskar-Morgenstern-Platz 1, A-1090 Vienna, Austria,
Vienna Research Platform on Accelerating
  Photoreaction Discovery, University of Vienna, W\"ahringerstra\ss e 17, 1090 Wien, Austria,
 \& Istituto di
  Matematica Applicata e Tecnologie Informatiche {\it E. Magenes}, via
  Ferrata 1, I-27100 Pavia, Italy
}
\email{ulisse.stefanelli@univie.ac.at}
\urladdr{http://www.mat.univie.ac.at/$\sim$stefanelli}

%\date{\today}


 \begin{abstract}
Curves of maximal slope are a reference
gradient-evolution notion in metric spaces and arise as variational
formulation of a vast class of nonlinear diffusion equations. Existence
theories for
curves of maximal slope are often based on minimizing-movements
schemes, most notably on the Euler scheme.
 We present here an alternative
minimizing-movements approach, yielding more regular discretizations,
serving as a-posteriori convergence estimator, and
allowing for  a simple convergence proof.  
 \end{abstract}

\maketitle

%%%%%%%%%%%%%%%%%%%%%%%%%%%%%%%%%%%%%%
\section{Introduction} 

Gradient-flow evolution in metric spaces has been the subject of
intense research in
the last years. Starting from the pioneering
remarks in \cite{DeGiorgi80}, the theory has been boosted
by the monograph by Ambrosio, Gigli, \& Savar\'e
\cite{Ambrosio08} and now encompasses existence and approximation results, as
well as long-time behavior, decay to equilibrium, and regularity \cite{Santambrogio}.

The applicative interest in evolution equations in metric spaces has been
revived by the seminal observations in \cite{JKO} and the work by Otto
\cite{Otto01} that a remarkably large class
of diffusion equations can be variationally reinterpreted as gradient
flows in Wasserstein spaces. More precisely, consider the nonlinear diffusion
equation
\begin{equation}
  \label{eq:introd2}
  \partial_t \rho - {\rm div}\big(  \rho\nabla (  V +
     F'(\rho) + W \ast \rho)   \big) = 0 \quad
  \text{in}\  \Rz^d \times (0,T).
\end{equation}
% \begin{equation}
%   \label{eq:introd}
%   \partial_t \rho -{\rm div}\, \left(\rho \nabla \left(\frac{\partial
%         \mathcal L}{\partial \rho} \right) \right)= 0 \quad \text{in} \ \ \Rz^d
%   \times (0,T) 
% \end{equation}
% where $\partial \mathcal L/\partial \rho\to \Rz$ is the variation of the functional
% $$\mathcal L (\rho) = \int_{\Rz^d}L(x,\rho(x),\nabla \rho(x))\,\dx.$$ 
Here,  $\rho=\rho(x,t)\geq 0$ is a time-dependent   density with fixed total mass
$\int_{\Rz^d}\rho(x,t)\, \dx=1$ and finite second moment 
$\int_{\Rz^d}|x|^2\rho(x,t)\, \dx<\infty$. Finally, 
$V : \Rz^d \to \Rz$ is a {\it confinement} potential, $F:[0,\infty) \to
\Rz$ is an {\it internal-energy density}, $W:\Rz^d \to \Rz$ is an {\it
  interaction} potential, and $\ast$ stands for the standard convolution in $\Rz^d$. 

Equation \eqref{eq:introd2}
can be variationally reformulated in terms of the gradient flow in the metric space $(\mathcal P_2(\Rz^d), W_2)$ of probability measures
with finite second moment endowed with the $2$-Wasserstein distance $W_2$
of the functional $\phi$ defined as
\begin{equation}
  \phi(u)= \int_{\Rz^d}V(x) \, {\rm d} u(x) +
\int_{\Rz^d}F(\rho(x)) \,\dx + \frac12\int_{\Rz^d\times \Rz^d} W(x{-}y)
\, {\rm d} (u \otimes u)(x,y)\label{eq:f}
\end{equation}
 if $u = \rho {\mathcal L}^d$ and $\phi(u)=\infty$ if $u$ is not absolutely
continuous with respect to the Lebesgue measure ${\mathcal L}^d$ in
$\Rz^d$
see \cite{Ambrosio08} and Section \ref{sec:Wass}.

The reference notion of solution to gradient flows in metric spaces is
that of {\it curves of maximal slope} \cite{DeGiorgi80}, see Definition \ref{def:curve}
below. This is based on a specific reformulation of \eqref{eq:introd2}
in form of a single scalar relation, featuring specific scalar quantities
playing the role of the norm of time derivative of the trajectory and
of the gradient of the energy, in the spirit of \eqref{eq:metric} below. 
Existence and decay
to equilibrium of curves of maximal slope for $\phi$ 
in $(\mathcal P_2(\Rz^d), W_2)$ are available, see
\cite{Ambrosio08,CMV,CMV2}, for instance. % and \cite{Agueh2003,rsss}, for instance.

In this paper, we focus on a novel time-discretization scheme for
gradient flows in metric spaces, falling within the class of {\it
  Minimizing Movements} in the sense of De Giorgi
\cite{Ambrosio95,DeGiorgi93}. Our theory is framed in abstract metric spaces, see Sections
\ref{sec:prem}-\ref{sec:p2}, and applied in linear and Wasserstein spaces in
Sections \ref{sec:appl0} and \ref{sec:Wass}, respectively. To keep this introductory
discussion as simple as possible, we present here the idea in the case of the doubly
nonlinear ODE system driven by a smooth potential $\phi$ on  $\Rz^d$, namely 
\begin{equation}
  |u'|^{p-2}u'+ \nabla \phi(u)=0 \quad \text{in} \  \times
(0,T)\label{ode}
\end{equation}
 for $p>1$, where the prime denotes time differentiation. 
This equation can be equivalently rewritten as
\begin{equation}
\phi(u(t)) + \frac1p\int_0^t |u'(r)|^p\dr +\frac1q\int_0^t 
|\nabla\phi(u(r))|^q\dr - \phi(u(0))=0 \quad \forall t \in (0,T)
\label{eq:metric}
\end{equation}
where now $q=p/(p-1)$ is 
conjugate to $p$. 
Note that the left-hand side above is always nonnegative, so that
\eqref{eq:metric} corresponds indeed to a so-called {\it null-minimization} principle: the left-hand side is minimized and
one checks that the minimum value is $0$. This approach has been lately
referred to as {\it De Giorgi's Energy-Dissipation} principle and has already
been applied in a variety of different contexts, including generalized
gradient flows \cite{Bacho,Rossi08}, rate-independent 
\cite{Mielke12,Roche} and GENERIC systems \cite{Peletier,generic_euler}, and optimal
control \cite{portinale}.

We  complement equation \eqref{ode} by specifying the initial
condition $u(0)=u^0$ for some $u^0 \in \Rz^d$. By introducing a time partition of $(0,T)$ with
uniform steps $\tau=T/N>0$, $N\in \Nz$ (note however that we consider  
nonuniform partitions below), and letting $ u_0=u^0$, the new minimizing-movements scheme reads
\begin{equation}
  \label{eq:2}
   u_i \in \argmin_u \left(\phi(u) +
    \frac{\tau^{1-p}}{p } \left|{u-u_{i-1}} \right|^p +\frac{\tau}{q}
    |\nabla \phi(u)|^q - \phi(u_{i-1}) \right)
\end{equation}
for $i=1,\dots,N$. With respect to the classical implicit Euler
method, scheme~\eqref{eq:2} includes an extra term featuring the
norm of the gradient. This modification with respect to Euler makes the function to be
minimized in \eqref{eq:2} a discrete and localized version of the
left-hand side in \eqref{eq:metric}. As such, scheme~\eqref{eq:2}
is nothing by the canonical {\it variational integrator} scheme
\cite{Hairer} associated
with the De Giorgi's Energy-Dissipation principle.

Compared to Euler,  the new
minimizing-movements scheme \eqref{eq:2} shows some distinguishing features. First of all, the direct occurrence of the gradient in \eqref{eq:2} entails additional
regularity of discrete solutions,  see
\eqref{eq:exreg}. As a matter of illustration, in the case of the {\it
  linear} heat equation ($p=2$) with homogeneous Dirichlet boundary conditions
scheme 
\eqref{eq:2} corresponds to solving the problem
$$\frac{u_i-u_i}{\tau}-\Delta u_i +\tau \Delta^2 u_i =0,$$
which is reminiscent of a {\it singular perturbation} of the Euler
scheme, see Section \ref{sec:illu}.

 Secondly, the exact correspondence of
\eqref{eq:2} to the left-hand side of \eqref{eq:metric} allows to check convergence of discrete solutions
without the need of introducing the so-called {\it De Giorgi's
  variational interpolation} function
\cite[Def. 3.2.1]{Ambrosio08}. 

Thirdly, in using a time discretization
to detect a minimum point of $\phi$ by iterating on the time steps,
the new scheme shows enhanced performance with respect to Euler for
large time steps, see \cite{generic_euler} and \eqref{prox} below. 

Finally, the functional under
minimization in \eqref{eq:2} may serve as an a-posteriori estimator for
the convergence of {\it any} discrete solution, regardless of the specific
method used to obtain it. In particular, one can resort to approximate
minimizers instead of true minimizers. 

The minimizing-movements scheme \eqref{eq:2} was already analyzed
in \cite{generic_euler} in the case of gradient flows in Hilbert
spaces. In particular,  convergence of the scheme for $\phi$ being a
$C^{1,\alpha}$ perturbation of a convex function and sharp, order-one error
estimates in finite dimensions can be found there. The case of curves of maximal slope in
metric spaces is also mentioned in \cite{generic_euler}, where
nevertheless the analysis is limited to $p=2$ and geodesically convex
potentials.

In this note, we extend the analysis of \cite{generic_euler} to the
case $p>1$ and to potentials $\phi$ being $(\lambda,p)$-generalized-geodesically convex for
$\lambda\in \Rz$.  More precisely, the combination of our main results, Theorems
\ref{thm:main1}-\ref{thm:main2}, entails that solutions to the new
minimizing-movements scheme \eqref{eq:2} in metric spaces, see
\eqref{eq:min}, converge to curves of maximal slope for all $p>1$, if
$\lambda \geq 0$, and for $p>2$, if $\lambda <0$.

In addition, in Theorem \ref{thm:main3} we are able to provide a
convergence result for not geodesically  convex functionals,
provided that some weak differentiability of its slope in form of a
generalized one-sided Taylor expansion condition holds, see~\eqref{eq:taylor}.

Before closing this introduction let us mention that alternative
time-discrete scheme with respect to Euler are available, also in the
nonlinear setting of metric spaces
\cite{Clement2,Matthes,Turinici,Tribuzio}.  We postpone an account on
the literature to Subsection
\ref{sec:literature}, for some preliminary material is needed to put
these contributions in perspective.

This is the plan of the paper. We introduce some notation and preliminaries in Section \ref{sec:prem}
and present our main convergence results in Section \ref{sec:main}. In
particular, assumptions are collected in Subsection
\ref{sec:assumptions} and statements are given in Subsection
\ref{sec:conv}. Some illustration of the theory on two linear
equations, both in finite and infinite dimensions, is in Subsection \ref{sec:illu}. The
convergence results are then proved in Sections
\ref{sec:p1}-\ref{sec:p3}. Eventually, we comment on the application
of the abstract theory in linear spaces in Section \ref{sec:appl0} and
in Wasserstein spaces  in Section \ref{sec:Wass}.

%%%%%%%%%%%%%%%%%%%%%%%%%%%%%%%%%%%%%%
\section{Preliminaries}\label{sec:prem}
\setcounter{equation}{0}

We briefly collect here some classical notation and preliminaries on evolution in metric
spaces, for completeness. The reader familiar with the classical
reference \cite{Ambrosio08} may consider moving directly to Section
\ref{sec:main}. 

In all of the following, $(U,d)$ denotes a complete metric
space and $\phi:U \to (-\infty,\infty]$ is a proper functional, i.e.,
the {\it effective domain} $D(\phi):=\{u \in U \ : \ \phi(u)
<\infty \}$ is assumed to be nonempty.

Let $p,\,q>1$ be given with $1/p+1/q=1$. A curve $u: [0,T]\to U$ is
said to belong to
$AC^p([0,T];U)$ if there exists $m\in
L^p(0,T)$ with 
\begin{equation}\label{metric_dev}
d(u(s),u(t))\leq \int_s^t m(r)\,\dr \quad \text{for all \ $0 \leq
s\leq t < T.$}
\end{equation}
If $u
\in AC^p([0,T];U) $, the limit
$$|u'|(t) := \lim_{s\to t}\frac{d(u(s),u(t))}{|t-s|}$$
exists for almost everywhere $t\in (0,T)$, see
\cite[Thm. 1.1.2]{Ambrosio08}, and is referred to as {\it metric
derivative} of $u$ at $t$. Moreover, the map
 $t \mapsto|u'|(t)$  is in $ L^p(0,T) $
and is minimal within the class of functions $m\in L^p(0,T)$
fulfilling \eqref{metric_dev}.

The {\it local slope} \cite{Ambrosio08,Cheeger99,DeGiorgi80} of $ \phi
$ at $ u \in D(\phi) $ is defined via
$$
\sl(u) := \limsup_{v \to u} \frac{(\phi(u) - \phi(v))^+}{d(u,v)}.$$ 
If $U$ is a Banach space and $\phi$ is Fr\'echet differentiable, we have that $\sl(u) = \|
{\rm D} \phi (u)\|_*$ (dual norm).

In the following, we will make use of the notion of {\it geodesic convexity} for $\phi$. More precisely, we call {\it (constant-speed) geodesic} any
curve $\gamma : [0,1]\to U$ such that $d(\gamma(t),\gamma(s)) =
(t-s)d(\gamma(0),\gamma(1))$ for all $0\leq s \leq t \leq T$ and 
we say that $\phi$ is {\it $(\kappa,p)$-geodesically convex} for 
  $\kappa\in \Rz$ if for all $ v_0,\, v_1 \in D(\phi)$ there exists a
  geodesic with
  $\gamma(0)=v_0 $ and $\gamma(1)=v_1$ such that 
\begin{align}
\phi(\gamma(\theta)) \leq \theta \phi(v_1) + (1-\theta) \phi(v_0) 
-\frac{\kappa}{p}\theta(1-\theta) d^p(v_0,v_1) \  \ \forall
  \theta \in [0,1]\label{def:l-p-geod-convex}
\end{align}
The definition is classical for $p=2$. For this $p$-extension see
\cite[Remark.~2.4.7]{Ambrosio08} or~\cite{Agueh2003}. Note that
geodesic convexity in particular implies that $U$ is a {\it geodesic
  space}, for each pair $v_0$, $v_1$ is connected by a geodesic. More
generally, we say
that $\phi$ is {\it $(\kappa,p)$-generalized-geodesically convex} if
\eqref{def:l-p-geod-convex} holds for some curve $\gamma$ connecting
$v_0$ and $v_1$, not necessarily being a geodesic. In this case, $U$
is implicitly assumed to be path-connected.

%  If $U$ is a Banach
% space, $\phi$ is $p$-geodesically convex iff it is $p$-convex in the
% following sense
% $$\phi(\theta v_1 + (1-\theta)v_0) \leq \theta \phi(v_1) + (1-\theta)
% \phi(v_0) - \frac{\kappa}{p}\theta(1-\theta) \| v_1- v_0\|^p$$
% for all $v_0,\,v_1
% \in D(\phi)$ and all $ \theta\in [0,1]$.

From \cite[Prop. 2.7]{rsss} we have that if $\phi$ is $(\kappa,p)$-geodesically convex and 
$d$-lower semicontinuous, the local slope $\sl$ is $d$-lower
semicontinuous as well. In addition, $\sl$ admits
the representation
\begin{equation}
 \sl(u)= \sup_{v \neq u} \left( \frac{\phi(u)
- \phi(v)}{d(u,v)} + \frac{\kappa}p d^{p-1}(u,v)\right)^+\quad
\forall u \in D(\phi).\label{repr-loc-slope}
\end{equation}
We denote by $D(\sl)$ the effective domain of $\sl$, namely,
$D(\sl)=\{u \in D(\phi) \ : \ \sl(u)<\infty\}$. Under the
above-mentioned geodesic convexity assumption, the local slope $\sl$ is  a {\it strong upper gradient}
\cite[Def.~1.3.2]{Ambrosio08}. Namely, for all $u\in AC^p([0,T];U)$,
the map
$r \mapsto \sl(r)$ is Borel and 
$$|\phi(u(t)) - \phi(u(s))| \leq \int_s^t \sl (u(r))\,\md(r) \, \dr
\quad \forall 0\leq s \leq t \leq T.$$
Note that, if $r \mapsto \sl (u(r))\md(r) \in L^1(0,T) $ the
latter entails that $\phi \circ u  \in W^{1,1}(0,T) $ and $(\phi
\circ u)' = \sl (u) \md$ almost everywhere in $(0,T)$.

Along with the above provisions, we specify the notion of
gradient-driven evolution as
follows.

\begin{definition}[Curve of maximal slope]\label{def:curve}
  The trajectory $u\in AC^p([0,T];U)$ is said to be a \emph{curve of
    maximal slope}
  if $\phi\circ u \in W^{1,1}(0,T)$ and 
  \begin{equation}
    \label{eq:curve}
    \phi(u(t)) + \frac1p\int_0^t \md^p(r) \, \dr  + \frac1q \int_0^t
    \sl^q(u(r))\, \dr =
    \phi(u(0)) \quad \forall t \in [0,T].
  \end{equation}
\end{definition}

%%%%%%%%%%%%%%%%%%%%%%%%%%%%%%%%%%%%%%%%%%%%%
\section{Main results}\label{sec:main}
\setcounter{equation}{0}
 
To each time partition $0=t_0<t_1<\dots< t_{N}=T$ we associate the
time steps $\tau_i = t_i - t_{i-1}$ and the diameter $\tau = \max
\tau_i$. Given the vector $\{u_i\}_{i=0}^N  \in U^{N+1}$ we define its
backward piecewise
constant interpolant $\ove u:[0,T] \to U$ on the time partition to be
\begin{align*}
&\ove u(0)= u_0\quad \text{and} \quad \ove u(t)=u_i \quad
\forall t \in (t_{i-1},t_i], \ i =1, \dots, N.
\end{align*}
Moreover, we define the piecewise constant function $|\haz
u'|:[0,T]\setminus\{t_0,\dots, t_{N}\} \to [0,\infty)$ as 
$$|\haz u'|(t) := \frac{d(u_{i-1},u_i)}{\tau_i} \quad \forall t \in
(t_{i-1},t_i), \ i =1, \dots, N.$$
The notation $|\haz u'|(t)$ alludes to the fact that in the Hilbert-space case the
latter is nothing but the norm of the time derivative of the piecewise
affine interpolant of the values $\{u_i\}_{i=0}^N$ on the time partition.

Our new minimizing-movements scheme is specified by means of the {\it incremental
functional} $G: (0,\infty) \times D(\phi) \times D(\sl) $ given by 
\begin{equation}
  \label{eq:G}
  \boxed{G(\tau,v,u):= \phi(u) + \frac{\tau^{1-p}}{p}d^p(v,u) +
  \frac{\tau}{q}\sl^q(u)- \phi(v).}
\end{equation}

In the setting of the assumptions specified later in Subsection
\ref{sec:assumptions},  for all $(\tau,v) \in (0,\infty) \times D(\phi)
$ the functional  $u\in D(\sl)\mapsto G(\tau,v,u)$ admits a
minimizer, possibly being not unique.  
We indicate the set of such minimizers by $M_G(\tau,v)$ and the minimum
value of $ G(\tau,v,\cdot)$  by $\haz G(\tau,v)$, namely, 
$$M_G(\tau,v):=\argmin_{u \in D(\sl)} G(\tau,v,u), \quad \haz G(\tau,v) :=
\min_{u \in D(\sl)}G(\tau,v,u).$$
With this notation, the new minimizing-movements scheme reads
\begin{equation}
 \boxed{u_0=u^0 \ \ \text{and} \ \ u_i \in M_G(\tau_i,u_{i-1})\quad \text{for}  \ i
=1,\dots,N,}\label{eq:min}
\end{equation}
for some given initial datum $u^0\in D(\phi)$.

For later purposes, we introduce also the  incremental
functional $E: (0,\infty) \times D(\phi) \times D(\phi)$  associated
with the classical backward Euler method
\begin{equation}
  \label{eq:F}
  E(\tau,v,u):= \phi(u) + \frac{\tau^{1-p}}{p}d^p(v,u) - \phi(v),
\end{equation}
as well as the corresponding notation
$$M_E(\tau,v):=\argmin_{u \in D(\phi)} E(\tau,v,u), \quad \haz E(\tau,v) :=
\min_{u \in D(\phi)} E(\tau,v,u).$$
In particular, the Euler method corresponds to the incremental problem
\begin{equation}
 u_0=u^0 \ \ \text{and} \ \ u_i \in M_E(\tau_i,u_{i-1})\quad \text{for}  \ i
=1,\dots,N.\label{eq:mine}
\end{equation}
In the context of Wasserstein spaces, see Section \ref{sec:Wass}, the latter is often referred to as
Jordan-Kinderlehrer-Otto scheme \cite{JKO}.

\subsection{Assumptions}\label{sec:assumptions}
 In this subsection, we fix our assumptions
and collect some comment. We start by asking that
\begin{equation}
  \label{eq:X}
  (U,d) \ \ \text{is a complete metric space}.
\end{equation}
In addition to the metric topology, $(U,d)$ is assumed to be endowed
with 
\begin{equation}
  \label{eq:sigma}
\text{a Hausdorff topology} \ \sigma, \  \text{compatible with the
  metric $d$}.
\end{equation}
The latter compatibility is intended in the following sense
\begin{equation}
  \label{eq:compat}
  u_n \sto u, \ v_n\sto v \  \ \Rightarrow \ \ d(u,v) \leq \liminf_{n
    \to \infty}  d(u_n,v_n)
\end{equation} 
and, in essence, means that $\sigma$ is weaker than the topology
induced by $d$. An early example for $\sigma$ complying with
\eqref{eq:sigma} is the topology  
induced by $d$. 
In applications it may however be useful to keep the two topologies
separate. In particular, if $U$ is a Banach space $\sigma$ is often
chosen to be some weak topology whereas $d$ usually corresponds to the strong
one. 

The initial datum is assumed to satisfy
\begin{equation}
  \label{eq:initial}
  u^0 \in D(\phi).
\end{equation}

We assume the proper potential $\phi:U \to (-\infty,\infty]$ to
be such that 
\begin{equation}
  \label{eq:phicomp}
\text{the sublevels of $\phi$ are sequentially
    $\sigma$-compact}.
\end{equation} 
The latter in particular entails that $\phi$ is sequentially
$\sigma$-lower semicontinuous and bounded from below. In the
following, we hence assume  
with no loss of generality that $\phi$ is nonnegative. Note however
that assumption \eqref{eq:phicomp} could be weakened by asking
compactness on $d$-bounded sublevels of $\phi$ only.

In addition, we ask
that 
\begin{align}
  &\sl \ \text{is a strong upper gradient for $\phi$ and it is sequentially}\nonumber\\
& \quad \text{$\sigma$-lower semicontinuous on $d$-bounded sublevels of $\phi$.}
  \label{eq:phisl}
\end{align}
The latter assumption could be weakened by developing the theory for some 
relaxation of $\sl$. Still,  \cite[Prop.~2.7]{rsss} ensures that
\eqref{eq:phisl} hold, as soon as $\phi$ is $(\lambda,p)$-geodesically
convex and $\sigma$ is the metric topology induced by $d$. 

In the setting of assumptions \eqref{eq:X}-\eqref{eq:phisl}, the 
solvability of the incremental minimization problem \eqref{eq:min}
follows from the Direct Method. Indeed, for all $\tau>0$ and $v\in
D(\phi)$ the incremental
functional $u \in D(\sl) \mapsto G(\tau,v,u)$ is coercive and 
lower semicontinuous by \eqref{eq:phicomp}-\eqref{eq:phisl}. We will later
check in \eqref{eq:slopeestimate2} that indeed 
\begin{equation} u\in M_G(\tau,v) \ \ \Rightarrow \ \ |\partial (\phi +\tau
\sl^q/q)|(u)<\infty. \label{eq:exreg}
\end{equation}
In particular, minimizers of $ G(\tau,v,\cdot)$ show additional
regularity. This extra regularity may be not preserved by the time-continuous limit.

Under the sole \eqref{eq:phicomp} the incremental Euler
minimization problem \eqref{eq:mine} is solvable as well. In
particular, for all $\tau>0$ and $v\in D(\phi)$ the functional $u \in D(\phi)
\mapsto E(\tau,v,u)$ admits a minimizer. 

Along the analysis, we will make reference to specific generalized geodesically
convex cases. In particular,
we may ask   for 
\begin{align}
&\exists \tau_*>0, \ \lambda \in \Rz \ \ \text{such that} \ \ \forall
  \tau \in (0,\tau_*), \ \forall v \in D(\phi)\nonumber\\
%&\inf_{v \in D(\phi)} F(\tau,u,v)>-\infty \ \ \text{and} \label{eq:F1}\\
&u \mapsto E(\tau,v,u) \ \ \text{is
                                                                            $(\kappa,p)$-generalized-geodesically convex}\nonumber\\
  & \text{with
  $\kappa=(p-1)\tau^{1-p}+\lambda$}.
\label{eq:F2}
\end{align}
Note that \eqref{eq:F2} holds if $\phi$ is $(\lambda,p)$-geodesically
convex and the  $p$-power of the distance is
$(p-1,p)$-geodesically convex. In case $p=2$, the $(1,2)$-geodesic
convexity of $u\mapsto d^2(u,v)/2$ qualifies nonpositively curved spaces in the Alexsandrov sense \cite{Alexandrov,Jost}. In
particular,  Euclidean and
Hilbert spaces, as well as Riemannian manifolds of
nonpositive sectional curvature \cite[Rem. 4.0.2]{Ambrosio08} fall
into this class. 

Condition \eqref{eq:F2} is more demanding for $p\not =2$. In fact, by
letting $\tau \to 0$ it implies
that the $p$-power of the distance is $(p-1,p)$-geodesically
convex. This is actually not the case in linear spaces, as one can
check already in $\Rz$, but see also \cite[Lem. 3.1]{AguehE}. Indeed, let $\theta=1/2$ and $v_0=-1$, $v_1=1$,
$\theta=1/2$ for $p>2$ and $v_0=0$, $v_1=1$ for $p<2$ in order to get
$$ \frac1p |\theta v_1 + (1-\theta) v_0|^p >\frac{\theta}{p}|v_1|^p + \frac{1-\theta}{p}|v_0|^p -
\theta(1-\theta)\frac{p-1}{p} |v_1-v_0|^p$$
contradicting $(p-1,p)$-geodesic convexity.
See 
\cite[Ex.~1, p.~55]{Jost} for some similar argument, proving 
the failure of $(1,2)$-geodesic
convexity of $(x_1,x_1) \in \Rz^2\mapsto (x_1^p + x_2^p)^{1/p}$.
In fact, condition \eqref{eq:F2} for $p\not =2$ is actually meaningful only in
spaces of qualified negative curvature. This is not the case for the
Wasserstein space $(\mathcal P_2(\Rz^d) , W_2)$, which is actually of
positive curvature, see Section \ref{sec:Wass}.
As we deal in 
Sections \ref{sec:appl0}-\ref{sec:Wass} with applications in linear
and Wasserstein spaces, condition~\eqref{eq:F2} is used there only for $p=2$.

In case of not geodesically convex potentials, we are still in the
position of providing a convergence result under the following
generalized one-sided Taylor-expansion condition on $\sl$
 \begin{align}
   &\exists \tau_*>0,  \, \forall C>0, \, \exists g:(0,\tau_*) \to [0,\infty] \ \text{with}  \
\frac{1}{\tau}\int_0^\tau g(r)\, {\rm d} r \searrow 0 \  \text{as} 
\ \tau \to 0  \ \text{such that}\nonumber\\
&\forall \tau \in (0,\tau_*), \ \forall v \in D(\sl)  \ \text{with} \ 
\max\{\phi(v),\tau \sl^q(v)\}\leq C, \ \forall u \in M_G(\tau,v) \nonumber\\[2mm]
   &                                                 \text{we have
     that} \ \ \sl^q(u) - |\partial (\phi +\tau\sl^q/q)|^q(u)  \leq
    g(\tau). \label{eq:taylor}
  \end{align}
Notice that the last inequality makes sense, for we have the
additional regularity
\eqref{eq:exreg}. We discuss some applications fulfilling condition
\eqref{eq:taylor} in Sections \ref{sec:appl0} and \ref{sec:Wass}.

A caveat on notation: In the following we use the same symbol $C$ in
order to indicate a generic positive constant, possibly depending on
data and changing from line to line. Where needed, dependencies are
indicated by subscripts.

%%%%%%%%%%%%%%%%%%%%%%%%%%%%%%%%%%%%%%%%%%%%%%%%%
\subsection{Convergence results}\label{sec:conv}
We are now ready to state our main results.

\begin{theorem}[Conditional convergence]\label{thm:main1} Under \eqref{eq:X}-\eqref{eq:phisl} let
  $\{0=t_0^n<t_1^n<\dots<t^{n}_{N^n}=T\}$ be a sequence of partitions
  with $\tau^n:=\max (t^n_i-t^n_{i-1}) \to 0$ as $n\to \infty$. Moreover, let
  $\{u_i^n\}_{i=0}^{N^n}$ be such that $u_0^n$ are $d$-bounded, $u_0^n \stackrel{\sigma}{\to}
  u^0$, $\phi(u_0^n)
  \to \phi(u^0)$, and
\begin{equation}
\sum_{i=1}^{N^n} (G(\tau_i^n,u_{i-1}^n,u_i^n))^+ \to 0 \ \
\text{as} \ \ n \to \infty.\label{eq:cond}
\end{equation}
Then, up to a not relabeled subsequence, we have that $\ove u^n(t) \sto u(t)$, where $u$ is a curve of
maximal slope with $u(0)=u^0$.
\end{theorem}

Note that the statement of Theorem
\ref{thm:main1}  does not require that $u^n_i\in
M_G(\tau^n_i,u^n_{i-1})$, namely that $\{u^n_i\}_{i=0}^{N^n}$ is a
solution of the new minimizing-movements scheme \eqref{eq:min}. In particular, Theorem
\ref{thm:main1} can serve as an a-posteriori tool to check the
convergence of time-discrete approximations, regardless of the method
used to generate them. In particular, the above conditional
convergence result directly applies to {\it approximate} minimizers,
namely solutions of
$$u_0^n=u^0 \quad \text{and} \quad G(\tau_i,u^n_{i-1},u^n_i) \leq
\inf G(\tau_i^n,u^n_{i-1},\cdot) +g^n_i\quad \text{for} \ i=1,\dots,N^n$$
(compare with \eqref{eq:min}) as long as $\sum_{i=1}^{N^n}g^n_i\to 0$
as $n\to \infty$. See \cite{Fleissner} for a result on approximate
minimizers of $E(\tau_i^n,u^n_{i-1},\cdot)$ instead.

The conditional convergence result of Theorem \ref{thm:main1} thus relies
on the possibility of solving the inequality $G(\tau_i^n , u_{i-1}^n,u_i^n)\leq
0$ up to a small, controllable error, and establishing some a priori
bounds on the discrete solution.
The validity of condition \eqref{eq:cond} is to be checked on the
specific problem at hand. In the specific case of $(\lambda,p)$-generalized-geodesically
convex functionals $\phi$ on a properly nonpositively curved space, condition
\eqref{eq:cond} actually holds for solutions of the new minimizing-movements scheme \eqref{eq:min}. This is the content of our second main result.

\begin{theorem}[Convergence in the geodesically convex case]\label{thm:main2}
  Under assump-\hfill \linebreak tions \eqref{eq:X}-\eqref{eq:phisl} and \eqref{eq:F2},  let
  $\{0=t_0^n<t_1^n<\dots<t^{n}_{N^n}=T\}$ be a sequence of partitions
  with $\tau^n:=\max (t^n_i-t^n_{i-1})<\tau_*$ and $\tau^n\to 0$ as
  $n\to \infty$. Moreover, assume that either $\lambda \geq 0$ or
  $p> 2$ in
  \eqref{eq:F2}. Then, solutions  $\{u_i^n\}_{i=0}^{N^n}$ of \eqref{eq:min} fulfill condition
\eqref{eq:cond}. Hence, $\ove u^n$ 
converges pointwise to a curve of maximal slope up to subsequences.
\end{theorem} 

% The restriction on $\lambda$ and $p$ in the statement of Theorem
% \ref{thm:main2} cannot be omitted. Indeed, in case $\lambda <0$ and
% $p<2$ the continuous problem may have no global solution (take
% $U=\Rz$, $\phi(u)=lambda u^2/2$ and $u^0>0$). Note nonetheless that
% the scheme may still show covnergence up to blow-up time.

We now turn to a convergence result in the not geodesically convex
case. Here, some stronger topological assumption, an approximation of
the initial datum, and the generalized one-sided
Taylor-expansion assumption \eqref{eq:taylor} for $\sl$ are necessary.

\begin{theorem}[Convergence without geodesic convexity]\label{thm:main3}
  Under assumptions \linebreak\eqref{eq:X}-\eqref{eq:phisl},  let $\sigma$ be the
  metric topology induced by $d$, $U$ be separable, and $\phi$ fulfill \eqref{eq:taylor}.
Moreover, let
  $\{0=t_0^n<t_1^n<\dots<t^{n}_{N^n}=T\}$ be a sequence of partitions
  with $\tau^n:=\max (t^n_i-t^n_{i-1})<\tau_*$, $(\tau_i^n - \tau_{i-1}^n)^+/\tau^n_{i-1} \leq
   C\tau^n$ for $i=2,\dots,N^n$,  and $\tau^n\to 0$ as
  $n\to \infty$. Choose $u^{0n} \in M_E(\tau^n,u^0)$. Then,
  solutions  $\{u_i^n\}_{i=0}^{N^n}$ of \eqref{eq:min} with $u_0^n = u^{0n}$  fulfill condition
\eqref{eq:cond}. Hence, $\ove u^n$
converges pointwise to a curve of maximal slope up to subsequences. 
\end{theorem}

Note that the one-sided nondegeneracy condition $(\tau_i^n - \tau_{i-1}^n)^+/\tau^n_{i-1} \leq
   C\tau^n$ in the statement is fulfilled if $i \mapsto
\tau^n_i$ in nonincreasing. In particular, it holds for uniform partitions.
In case $u^0\in D(\sl)$ no approximation of the initial datum as in
Theorem \ref{thm:main3} is actually needed.

Theorems \ref{thm:main1}, \ref{thm:main2}, and  \ref{thm:main3} are proved in Sections
\ref{sec:p1},  \ref{sec:p2}, and \ref{sec:p3}, respectively.

\subsection{An illustration on linear equations}\label{sec:illu}
The focus of our theory is on nonlinear problems. Still, as a way of
illustrating the results, we present here two linear ODE and PDE
examples. Nonlinear applications are then discussed in Sections
\ref{sec:appl0}-\ref{sec:Wass} below. 

Let us start from the finite-dimensional example of the gradient flow in
$(\Rz^d,|\cdot|)$ of
$\phi(u)=\lambda |u|^2/2$ with $\lambda \in \Rz$ and take $p=2$. In this case, the incremental functional $G$ reads
$$G(\tau,v,u) = \frac{\lambda}{2}|u|^2 + \frac{1}{2\tau}|u-v|^2 +
\frac{\tau \lambda^2}{2}|u|^2 - \frac{\lambda}{2}|v|^2.$$
For all $v\in \Rz^d$ given, the latter can be readily minimized,
giving the only minimum point $u = v/(1+\lambda \tau + \lambda^2 \tau^2)$. Correspondingly,
the minimal value $\haz G(t,v)$ can be checked to be 
\begin{equation}\haz G(t,v) ={}
-\frac{|v|^2\lambda^3\tau^2}{2(1+\lambda \tau + \lambda^2
  \tau^2)}.\label{eq:fin}
\end{equation}
If $\lambda \geq 0$ the minimal value is nonpositive and condition
\eqref{eq:cond} trivially holds. If $\lambda <0$, the minimal value scales as
$\tau^2$ and condition \eqref{eq:cond} still holds. Indeed,
by letting 
\begin{equation}
  \label{r}
  r^n :=\sum_{i=1}^{N^n} \big( G(\tau_i^n,u_{i-1}^n,u_i^n)\big)^+
\end{equation}
we have that  
\begin{equation} r^n =
 \sum_{i=1}^{N^n} \frac{|u_{i-1}^n|^2(\lambda^-)^3(\tau_i^n)^2}{2(1+\lambda \tau_i^n + \lambda^2
  (\tau_i^n)^2)}\leq C \max_i|u_i^n|^2 \tau^n\label{lin}
\end{equation}
where we tacitly assumed that $\lambda^-\tau^n\leq \lambda^-\tau_*<1$ and we used the standard notation for
the negative part $\lambda^- = \max\{0,-\lambda\}$. Condition
\eqref{eq:cond} hence follows as soon as $ \max_i|u_i^n|$ stays
bounded with respect to $n$, which happens to be the case as the
evolution takes place in the finite time interval $[0,T]$.

\begin{figure}[ht]
\centering
\pgfdeclareimage[width=65mm]{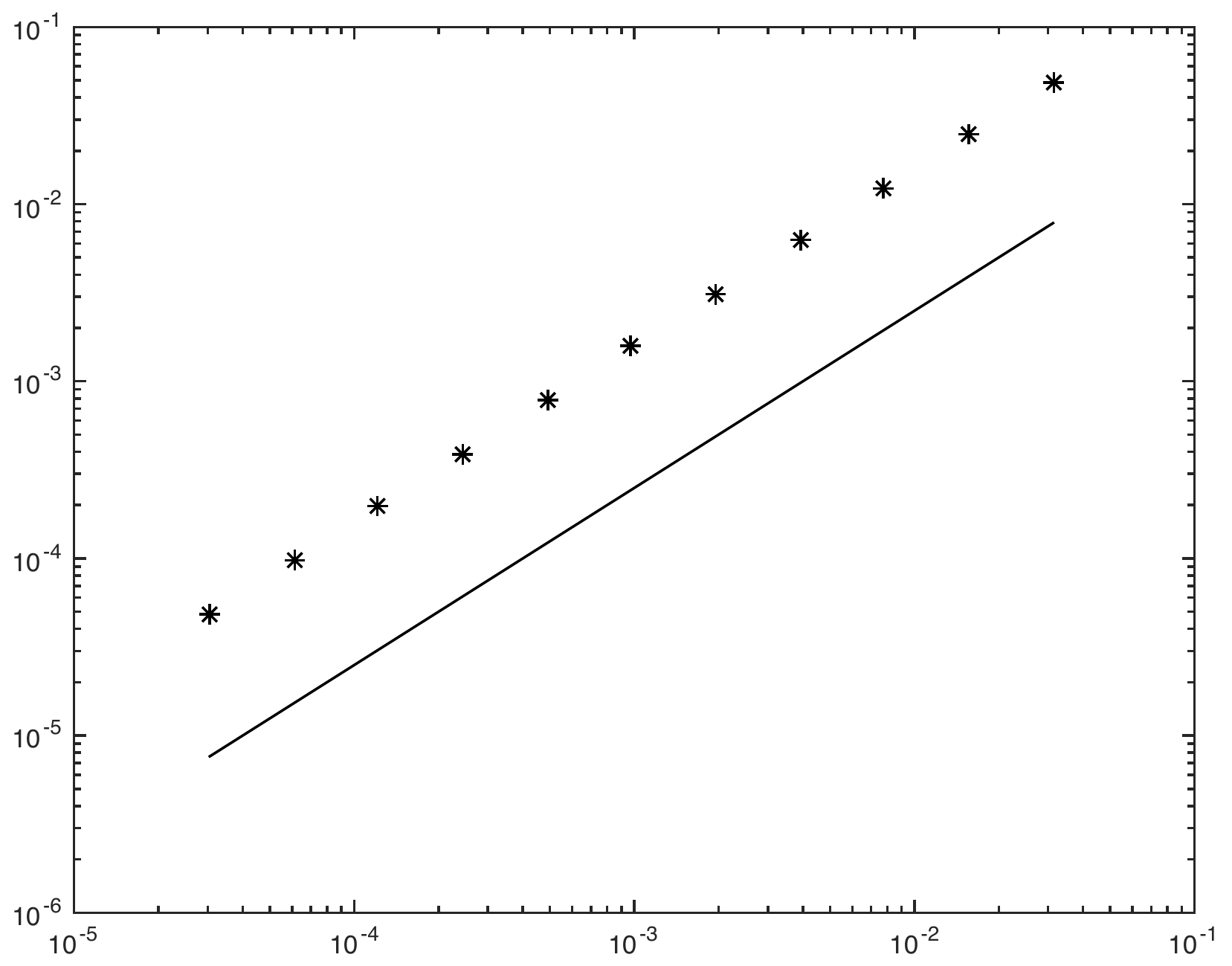}{r}  
\pgfuseimage{r}%\pgfuseimage{fig2}
\caption{Values $r^n$ from \eqref{r} against $\tau^n$ (stars)
  with respect to order $1$ (solid) in log-log scale.}
\label{figure1}
\end{figure}

In fact, the order of convergence in \eqref{lin} is sharp, as
illustrated in  Figure
\ref{figure1} for the choice $d=1$, $\lambda=-1$,  $  u^0=1$, $T=1$. 
Here, $r_n$ in computed for  the uniform partition
$\tau^n_i=\tau^n=2^{-n}$, $n=1, \dots, 12$ or, equivalently, for $N^n=2^n$. %The value  $r_n$  is
%plotted against $\tau^n$ in a loglog scale.
 
% Convergence of the new minimizing-movement scheme for
% $\phi(u)=\lambda|u^2|/2$ is indeed guaranteed by Theorem
% \ref{thm:main2} in the convex case $\lambda\geq 0$ and by Theorem
% \ref{thm:main3} in the nonconvex case $\lambda<0$, for $\phi$ is
% smooth.

\begin{figure}[ht]
\centering
\pgfdeclareimage[width=65mm]{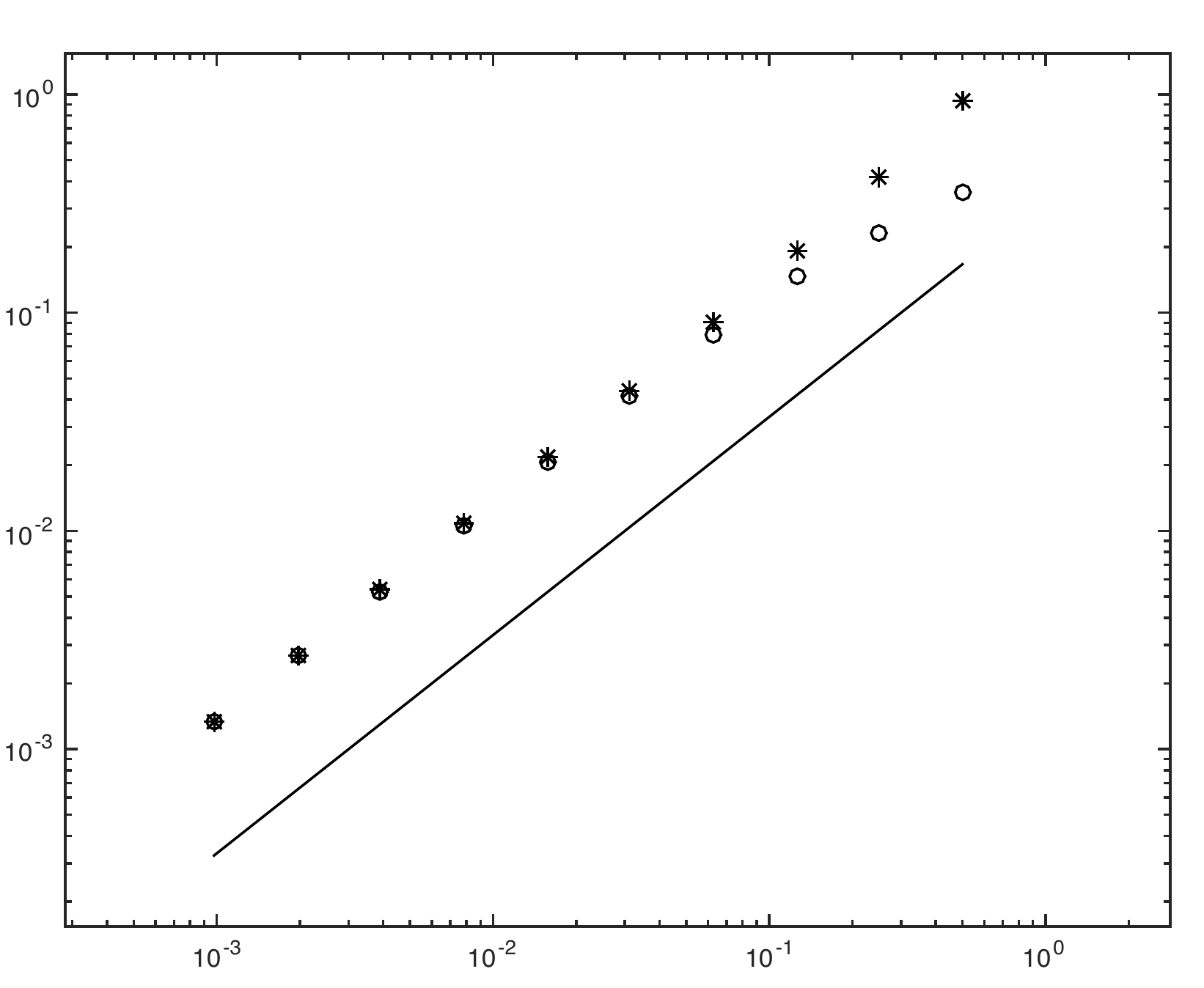}{comparison}  
\pgfuseimage{comparison}%\pgfuseimage{fig2}
\caption{$L^\infty$ error with respect to $\tau^n$ for the new
  minimizing-movements scheme (stars) and the Euler scheme (dots) in log-log
  scale. The solid line represents order $1$.}
\label{figure2}
\end{figure}

On a uniform partition of time step $\tau>0$, the solution of the new
minimizing movement scheme $\{u_i\}$ and the solution $\{u_i^e\}$ of the Euler scheme read
\begin{equation}
  u_i = \frac{u_0}{(1+\lambda \tau+ \lambda^2 \tau^2)^i}\quad \text{and} \quad u_i^e
= \frac{u_0}{(1+\lambda \tau)^i},\label{eq:veresol}
\end{equation}
respectively. It is hence a standard matter to compute
\begin{equation} |u_i - u^e_i| = |u_0|\left| \frac{(1+\lambda \tau)^i -
    (1+\lambda \tau+ \lambda^2 \tau^2)^i}{
    (1+\lambda \tau+ \lambda^2 \tau^2)^i (1+\lambda \tau)^i}\right
| \label{eq:standard}
\end{equation}
which scales like $\tau^2$ as $\tau \to 0$. As the Euler scheme is of first
order, the same holds true for the new minimizing-movements scheme, see
Figure \ref{figure2} for $d=1$, $\lambda=-1$, $u_0=1$. Indeed, Figure
\ref{figure2} shows that this order is sharp. Note in fact that  the new
minimizing-movements scheme is proved in
\cite[Prop. 4.3]{generic_euler} to be of first order for all
nonnegative potentials
$\phi$ in $C^2$ in finite dimensions.

Assume now to be interested in computing the minimum of $\phi$ by
following the discrete scheme for a fixed number $m$ of iterations, a
classical strategy in optimization \cite{Combettes,Rockafellar}. In
the specific case of our ODE example we compute from
\eqref{eq:veresol}  
\begin{equation} \phi(u_m) = \frac{\lambda}{2(1+\lambda\tau +\lambda^2
  \tau^2)^{2m}} \quad  \text{and} \quad \phi(u_m^e) =
\frac{\lambda}{2(1+\lambda\tau)^{2m}}.\label{prox}
\end{equation}
Due to the presence of the extra
term $\lambda^2 \tau^2$ in the denominator, the new
scheme is advantageous with respect to Euler as for reduction of the potential after a fixed number of
iterations. Note that this effect is enhanced by choosing {\it large}
time steps.

Let us move to an infinite-dimensional example by considering the standard heat
equation on the space time cylinder $\Omega \times (0,T)$ where
$\Omega \subset \Rz^d$ is a smooth, open, and  bounded set and homogeneous
Dirichlet conditions are imposed (other choices being of course
possible). We classically reformulate this as the gradient flow in
$(L^2(\Omega), \| \cdot \|)$, of the {\it Dirichlet} energy
$$\phi(u)
=
\left\{ 
\begin{array}{ll}\disp\frac12
\int_\Omega  |\nabla u(x)|^2 \, \dx \quad& \text{for} \ \ u
\in H^1_0(\Omega)\\
\infty&\text{elsewhere in} \  L^2(\Omega). 
\end{array}
\right.
$$
where $\| \cdot \|$ is the norm corresponding to the natural $L^2$ scalar
product $(\cdot, \cdot)$.  
In this case, we have that $\partial \phi (u) = - \Delta u$
with $D(\partial
\phi) =  H^2(\Omega) \cap H^1_0(\Omega)$. The symbol $\partial$
indicates the {\it subdifferential} in the sense of convex analysis   \cite{Brezis73}. In particular, $\partial
\phi$ is single-valued and $|\partial \phi|(u) = \| \Delta u \|$
for all $u\in D(\partial \phi)$. The incremental functional $G:
(0,\infty) \times  H^1_0(\Omega) \times H^2(\Omega) \cap H^1_0(\Omega)$ hence reads
$$G(\tau,v,u) = \int_\Omega \left( \frac12 |\nabla u|^2 +
  \frac{1}{2\tau}|u-v|^2+ \frac{\tau}{2} |\Delta u|^2 - \frac12
  |\nabla v|^2\right) \dx .$$
For all $v\in H^1_0(\Omega)$ given, the latter can be readily
minimized in $H^2(\Omega) \cap H^1_0(\Omega)$. Given linearity one can easily identify
the subgradient of $ u \mapsto G(\tau,v,u)$ as
$$(\partial G(\tau,v,\cdot))(u) = -\Delta u + \frac{u-v}{\tau} +\tau
\Delta^2 u   $$ and
$D(\partial G(\tau,v,\cdot)) = \{ u \in H^4(\Omega)\cap H^1_0(\Omega)
\, : \, \Delta u =0 \ \text{on} \ \partial \Omega\}$. Hence, 
  the minimizer $u$ of $G(\tau,v,\cdot)$ solves
  \begin{align}
    u -\tau\Delta u +\tau^2 \Delta^2 u =v \ \ \text{a.e. in} \ \
    \Omega,\quad 
    u=\Delta u = 0 \ \ \text{on} \ \ \partial \Omega. \label{eq:heat}
  \end{align}
  The latter is reminiscent of a singular perturbation of
 \begin{align}
    u^e -\tau\Delta u^e   =v \ \ \text{a.e. in} \ \
    \Omega,\quad 
    u^e = 0 \ \ \text{on} \ \ \partial \Omega,\label{eq:heat2}
 \end{align}
 corresponding instead to the incremental step of the Euler
 scheme.

 Let now $\{w^k\}$ be a complete orthonormal basis of $L^2$
of eigenfunctions of $-\Delta$ with homogeneous Dirichlet boundary
conditions, namely, $w^k\in H^2(\Omega) \cap H^1_0(\Omega)$ with $w^k
\not = 0$ and
 $-\Delta w^k = \lambda^k w^k$ for some $\lambda^k>0$. By inserting in
 \eqref{eq:heat}-\eqref{eq:heat2} $u
 = \sum_ku^kw^k$, $u^e
 = \sum_k(u^e)^kw^k$, and $v
 = \sum_kv^kw^k$ for $u^k := (u,w^k)$, $(u^e)^k := (u^e,w^k)$, and $v^k:=(v,w^k)$, respectively,
 we get that 
 $$ u^k = \frac{v^k}{1+ \tau \lambda^k + (\tau\lambda^k)^2} \quad
 \text{and}  \quad (u^e)^k = \frac{v^k}{1+ \tau \lambda^k}.$$
 In particular, by arguing as in \eqref{eq:fin} one readily checks that 
$$\haz G(\tau,v) = -\sum_k \frac{|v^k|^2(\lambda^k)^3\tau^2}{2(1+ \tau
  \lambda^k + (\tau\lambda^k)^2)} \leq 0$$
and condition \eqref{eq:cond} holds. By iterating on the time steps,
the solution $\{u_i\}$ of the new minimizing movement scheme and that
$\{u^e_i\}$ of the Euler scheme read $u_i = \sum_k u^k_i w^k$ and
$u^e_i = \sum_k (u^e)^k_i w^k$ where
$$ u^k_i = \frac{(u^0)^k}{(1+ \tau \lambda^k + (\tau\lambda^k)^2)^i} \quad
 \text{and}  \quad (u^e_i)^k = \frac{(u^0)^k}{(1+ \tau \lambda^k
   )^i}$$
 and $(u^0)^k:=(u^0,w^k)$. Proceeding as in \eqref{prox} one computes
 \begin{align*}
  & \phi(u_m) = \frac12 \sum_k \lambda^k (u^k_m)^2 = \frac12 \sum_k
    \frac{\lambda^k ((u^0)^k)^2}{(1+\tau \lambda^k + (\tau \lambda^k))^{2m}}, \\
   &\phi(u_m) = \frac12 \sum_k \lambda^k ((u^e_m)^k_m)^2 = \frac12 \sum_k
   \frac{\lambda^k ((u^0)^k)^2}{(1+\tau \lambda^k)^{2m}}
 \end{align*}
 and the same observations as in the ODE case on the effectiveness of the reduction of
 the potential for a fixed number of iterations apply.

\subsection{Literature}\label{sec:literature}
 
Before moving on, let us record here some other alternatives to the
Euler scheme, specifically focusing on the case 
$p=2$.

Legendre and Turinici advance in \cite{Turinici} the {\it midpoint} scheme
\begin{align*}
  u_i \in \argmin_{u} \Bigg( \inf \Bigg( 2\phi(w) + \frac{1}{2\tau}
  d(u,u_{i-1})\ : \ w \in \Gamma(u,u_{i-1}) \Bigg) \Bigg)
\end{align*}
where
$$\Gamma(u,u_{i-1}) = \{\gamma (1/2) \  :  \ \gamma:[0,1]\to U \
\text{geodesic with}   \ \gamma(0)=u_{i-1} \ \text{and} \
\gamma(1)=u\}.$$
By assuming \eqref{eq:phicomp}-\eqref{eq:phisl}, as well as some additional
closure property relating to the specific structure of the set $\Gamma$, they prove that this midpoint scheme is solvable
and convergent.

A variant of this scheme is also proposed in \cite{Turinici} in the specific case of nonbranching geodesic spaces,
namely, spaces where any two points are connected by a unique
geodesic. In these spaces, for all $w$ and $u_{i-1}$
there exists a unique $u$ such that $w\in \Gamma(u,u_{i-1})$. An {\it
  extrapolated} version of the Euler scheme is hence defined by the relations
$$ u^e_{1/2} \in \Gamma(u_i,u_{i-1}) \ \ \text{where} \ \ u^e_{1/2}
\in M_E(\tau/2,u_{i-1}).$$
Albeit not purely variational, this scheme is based on the solution of
the Euler scheme with halved time step.

Matthes and Plazotta \cite{Matthes} address a variational version
of the Backward Differentiation
Formula (BDF2) method, namely,
$$u_i \in \argmin_{u\in
  D(\phi)}\left(\frac{1}{\tau}d^2(u,u_{i-1})-\frac{1}{4\tau}
  d^2(u,u_{i-2})+\phi(u) \right) \quad \text{for} \ i=2,\dots,N$$
where now both $u_0$ and $u_1$ are given. Under some lower
semicontinuity and convexity conditions, it is proved in
\cite{Matthes} that the scheme admits a solution, whose
piecewise-in-time interpolant converges to a curve of maximal slope
with rate $\tau^{1/2}$. It also shown that under natural regularity
assumptions on the limiting time-continuous curve of maximal slope,
the convergence rate can be $\tau$ at best.

Perturbations of the Euler method of the form
$$u_i \in \argmin_{u\in
  D(\phi)}\left(\frac{a_i^\tau}{2\tau}d^2(u,u_{i-1}) +\phi(u) \right) \quad
\text{for} \ i=2,\dots,N,$$
are considered by Tribuzio in \cite{Tribuzio}. Here, one is given the
sequence of positive weights defined as  $a_i^\tau=a^\tau(i\tau)$ for some
functions $a^\tau:
(0,\infty) \to (0,\infty)
$. This generalization with respect to the classical Euler scheme yields a modification of
the metric as time evolves. By asking $1/a^\tau$ to be locally equiintegrable with respect to
$\tau$, one can prove that minimizers converge to curves of maximal
slope according to a specific time-dependent limiting
metric. Under some more general assumptions on $a^\tau$, discontinuous
evolutions can also be obtained. These can be proved to be capable
of exploring the different wells of a multiwell potential $\phi$.

Let us also mention the approach \`a la Crandall-Liggett by Cl\'ement
and Desch~\cite{Clement1,Clement2}, see also~\cite{Clement3}, who recursively define $u^n_0=u^0$ and $u^n_i =
J(u^n_{i-1})$ for $i=1,\dots,N = T/\tau$,
where   $J(u^n_{i-1})$  is the set of points $u\in D(\phi)$ fulfilling the
inequality
$$\frac{1}{2\tau} d^2(u,w) - \frac{1}{2\tau}d^2(u_{i-1}^n,w)
+\frac{1}{2\tau}d^2(u,u_{i-1}^n) + \phi(u) \leq \phi(w) \quad \forall
w \in D(\phi).$$
Such points exist for $\phi$ geodesically convex and the corresponding
interpolants $\ove u^n$ converge to {\it evolutionary variational
  inequality} solutions~\cite{Muratori}, a specific class of
curves of maximal
slope.

% Often gradient flows are employed to compute the minimum of
% functionals We can compare the performances of the new minimizing-movement
% scheme and the Euler scheme in reducing the potential $\phi(u)=u^2/2$
% $(d=1)$ along a fixed number of iterations as a function of the time
% step, see Figure \ref{figure3}.
% \begin{figure}[ht]
% \centering
% \pgfdeclareimage[width=65mm]{reduction}{reduction}  
% \pgfuseimage{reduction}%\pgfuseimage{fig2}
% \caption{Values of the potential $\phi(u)=u^2/2$
% $(d=1)$ after 20 iterations from $u_0=1$ with respect to $\tau^n$: the new
%   minimizing-movements scheme (solid) and the Euler scheme (dashed) in log-log
%   scale.}
% \label{figure2}
% \end{figure}

% Again in the scalar case $d=1$ one can repeat the experiment by
% choosing $\phi(u)=\lambda |u|^r/r$ for $r>2$. If
% $\lambda \geq 0$ the potential $\phi$ is convex and the corresponding
% $r^n$ is nonpositive. If $\lambda<0$, the potential $\phi$ is not
% $\lambda$-convex anymore and Figure \ref{figure1} right shows that $r^n$ does not
% converge to $0$ as $n\to \infty$. This proves that
% condition \eqref{eq:cond} can be violated if 
% $\lambda$-convexity is not assumed. In particular, $\lambda$-convexity
% is a necessary condition for the validity of the statement of Theorem \eqref{thm:main2}

% On the other hand, let us mention that the convergence of Figure
% \ref{figure1} left seems to hold also for $p\in (1,2]$. This suggests
% that the restriction $p>2$ in Theorem \ref{thm:main2} for the $\lambda$-convex case with $\lambda
% <0$ could be inessential.

%%%%%%%%%%%%%%%%%%%%%%%%%%%%%%%%%%%%%%
\section{Conditional convergence} \label{sec:p1}

This section is devoted to the proof of Theorem
\ref{thm:main1}. The ingredients of the argument are quite classical. Still, as already
mentioned, the current minimizing-movement setting of \eqref{eq:min} expedites the proof,
for there is no need to resort to the De~Giorgi variational
interpolant 
\cite[Def. 3.2.1]{Ambrosio08}.

Let $\{u_0^n\}$ be $d$-bounded  with $u_0^n \stackrel{\sigma}{\to} u^0$ and $\phi(u_0^n)
  \to \phi(u^0)$ fulfill \eqref{eq:cond}. We have that 
  \begin{align}
&\phi(\ove u^n(t_m^n)) + \frac1p\int_0^{t^n_m} |(\haz u^n)'| ^p(r)\, \dr +\frac1q
    \int_0^{t^n_m} \sl^q(\ove u^n(r))\, \dr \nonumber\\
&\quad = \phi(u_m^n) + \frac1p\sum_{i=1}^m (\tau_i^n)^{1-p}
    d^p(u_{i-1}^n,u_{i}^n) + \frac1q\sum_{i=1}^m \tau_i^n
  \sl^q(u_i^n)\nonumber\\
& \quad= \sum_{i=1}^m G(\tau_i^n,u_{i-1}^n,u_{i}^n) + \phi(u_0^n).\label{eq:pass}
  \end{align}
Condition \eqref{eq:cond} ensures that the above right-hand is bounded
independently of $m=1,\dots,N^n$ and $n$.
A first consequence of estimate \eqref{eq:pass} is that $\{u_m^n\}$ is
$d$-bounded  independently of $m=1,\dots,N^n$ and $n$. Indeed, one
has that
\begin{align*}
 &  d^p(u^n_0,u^n_m) \leq 2^{p-1 }\sum_{i=1}^m  d^p(u^n_{i-1},u^n_{i}) \leq
  2^{p-1} (\tau^n)^{p-1} \sum_{i=1}^m
  (\tau_i^n)^{1-p}d^p(u^n_{i-1},u^n_{i})\nonumber\\
  &\leq \quad 2^{p-1}(\tau^n)^{p-1}p \left(\sum_{i=1}^m G(\tau_i^n,u_{i-1}^n,u_{i}^n) + \phi(u_0^n)\right).
\end{align*}
The right-hand side is bounded 
independently of $m=1,\dots,N^n$ and $n$. Since $\{u_0^n\}$ are $d$-bounded, the $d$-boundedness of $\{u_m^n\}$ follows.

As the sublevels of $\phi$ are
sequentially $\sigma$-compact, one can apply the extended
Ascoli-Arzel\`a Theorem from \cite[Prop. 3.3.1]{Ambrosio08} and find a not
relabeled subsequence $\{\ove u^n\}$ such that $\ove
u^n\stackrel{\sigma}{\to} u$ pointwise, where $u:[0,T]
\to U$, and $|(\haz u^n)'| \to m$ weakly in $L^p(0,T)$. In particular, we
have that $u(0)=\lim_{n\to \infty} u^n(0) = \lim_{n\to
  \infty}u^n_0=u^0$. For
all $0 < s \leq t < T$, define $s^n=\max\{t^n_i \, : \, t_i^n<s\}$ and
$t^n=\min\{t^n_i\, : \, t < t^n_i\}$. Then,
$$d(u(s),u(t)) \stackrel{\eqref{eq:compat}}{\leq} \liminf_{n \to \infty} d(\ove u^n(s), \ove u^n(t)) \leq
\liminf_{n \to \infty} \int_{s^n}^{t^n} |(\haz u^n)'|(r)\, \dr = \int_s^t m(r) \,\dr.$$
This entails that $u \in AC^p([0,T];U)$ since we just checked that the
function $m
\in L^p(0,T)$ fulfills \eqref{metric_dev}.
As $|u'|$ is the minimal
function in $L^p(0,T)$ fulfilling \eqref{metric_dev}, we also 
have that  $|u|\leq m$ almost everywhere and
$$ \int_0^t |u'|^p(r)\,\dr  \leq \int_0^t m^p(r) \,\dr \leq  \liminf_{\tau
  \to 0} \int_0^t |(\haz u^n)'|^p(r)\,\dr 
\quad \forall t>0. $$
For all fixed $t\in (0,T]$, choose $t^n_m = t^n$ in \eqref{eq:pass} in order to get that 
 \begin{align*}
&\phi(\ove u^n(t)) + \frac1p\int_0^{\ove t^n(t)} |(\haz u^n)'| ^p(r)\, \dr +\frac1q
    \int_0^{\ove t^n(t)} \sl^q(\ove u^n(r))\, \dr \\
&\quad \stackrel{\eqref{eq:pass}}{\leq} \sum_{i=1}^{N^n}
  (G(\tau_i^n,u_{i-1}^n,u_{i}^n))^+ + \phi(u_0^n). 
  \end{align*}
Owing to the sequential $\sigma$-lower semicontinuity of $\phi$ and
$\sl$, see \eqref{eq:phicomp}-\eqref{eq:phisl}, we can pass to the $\liminf$ in
the latter and, using again condition \eqref{eq:cond} and the
fact that $\phi(u^n_0) \to \phi(u^0)$, we obtain 
\begin{equation}\phi(u(t)) + \frac1p \int_0^t \md^p(r)\, \dr + \frac1q \int_0^t
\sl^q(u(r))\, \dr \leq \phi(u(0)) \quad \forall t \in
[0,T].\label{eq:actually}
\end{equation}
As $\sl$ is a strong upper gradient for $\phi$ by \eqref{eq:phisl}, we
have that 
\begin{align*}
&\phi(u(0)) \leq \phi(u(t)) + \int_0^t \sl(u(r))\,\md(r)\, \dr \\
&\quad \leq \phi(u(t)) + \frac1p \int_0^t \md^p(r)\, \dr + \frac1q \int_0^t
\sl^q(u(r))\, \dr 
\end{align*}
so that \eqref{eq:actually} is actually an equality and $u$ is a curve of maximal slope
in the sense of Definition \ref{def:curve}.

%%%%%%%%%%%%%%%%%%%%%%%%%%%%%%%%%%%%%%
\section{Convergence in the  geodesically convex
  case} \label{sec:p2}
\setcounter{equation}{0}

We now turn to the proof of Theorem \ref{thm:main2}. 

Recall that for all $\tau^n_i>0$ and $v\in D(\phi)$ the functional $u \in D(\phi)
\mapsto E(\tau_i^n,v,u)$ admits a minimizer. 
We first prove a $p$-variant for $p>1$ of the slope estimate \cite[Lem. 3.1.3,
p. 61]{Ambrosio08}, which was originally proved for $p=2$. In particular, we
aim at the following
\begin{equation}
  \label{eq:slopeestimate}
  \sl(u) \leq (\tau_i^n)^{1-p} d^{p-1}(v,u) \quad \forall u \in M_E(\tau_i^n,v).
\end{equation}
Note that this estimate is already mentioned in \cite[Rem. 3.1.7]{Ambrosio08} without
proof. We give an argument here.
Let $w\in D(\phi)$ be given. From the
minimality $E(\tau_i^n,v,u) \leq E(\tau_i^n,v,w) $ we deduce that 
\begin{align*}
  &\phi(u) - \phi(w) \leq \frac{(\tau^n_i)^{1-p}}{p}\Big(d^p(v,w)
  -d^p(v,u) \Big) \\
&\quad\leq \frac{(\tau^n_i)^{1-p}}{p}\Big( \big( d(u,w)+d(v,u) \big)^p
  -d^p(v,u) \Big) \\
&\quad= \frac{(\tau^n_i)^{1-p}}{p} \left( 
\sum_{k=0}^\infty \binom{p}{k} d^k(u,w)\,d^{p-k}(v,u) -
  d^p(v,u)\right)\\
&\quad = d(u,w)  \frac{(\tau^n_i)^{1-p}}{p}  
\sum_{k=1}^\infty \binom{p}{k} d^{k-1}(u,w)\,d^{p-k}(v,u)
\end{align*}
where we have made use of the generalized binomial formula and the
generalized binomial coefficients
$$\binom{p}{k} = \frac{p(p-1)\dots (p-k+1)}{k!}.$$
Assume now that $w\not = u$, divide by $d(u,w)$, and compute the
$\limsup$ as $w \to u$ in order to get 
\begin{align*}
&\sl(u) = \limsup_{w \to u}\frac{\big(\phi(u) - \phi(w)\big)^+}{d(u,w)} \leq \limsup_{w \to u}\frac{(\tau^n_i)^{1-p}}{p}  
\sum_{k=1}^\infty \binom{p}{k} d^{k-1}(u,w)\,d^{p-k}(v,u) \\
&\quad =   \frac{(\tau^n_i)^{1-p}}{p}  
 \binom{p}{1}  d^{p-1}(v,u) = {(\tau^n_i)^{1-p}} 
   d^{p-1}(v,u)
\end{align*}
so that \eqref{eq:slopeestimate} holds. Above, we have used the fact
that 
\begin{align*}
  &0\leq \lim_{w\to u} \sum_{k=2}^\infty \binom{p}{k}
  d^{k-1}(u,w)\,d^{p-k}(v,u) \leq \lim_{w\to u}d(u,w) \sum_{k=2}^\infty
  \binom{p}{k}d^{p-k}(v,u)\\
&\quad = \lim_{w\to u} d(u,w) \left( (1+d(v,u))^p - \binom{p}{1}d^{p-1}(v,u) - \binom{p}{0}d^p(v,u)\right) = 0.
\end{align*}

Let now $u^e \in D(\phi)$ be a minimizer of $u\mapsto
E(\tau_i^n,u_{i-1}^n,u)$. Taking into account the 
convexity assumption \eqref{eq:F2}, let $\gamma:[0,1]\to U$ be a
curve with $\gamma(0)=u_{i-1}^n$ and $\gamma(1) =
u^e$, so that  
\begin{align*}
  &E(\tau_i^n,u_{i-1}^n,u^e) \leq  E(\tau_i^n,u_{i-1}^n,\gamma(\theta)) \\
&\quad \stackrel{\eqref{eq:F2}}{\leq}
  \theta  E(\tau_i^n,u_{i-1}^n,u^e) + (1-\theta)  E(\tau_i^n,u_{i-1}^n,
  u_{i-1}^n) \\
&\quad- \theta (1-\theta)\frac{(p-1)(\tau_i^n)^{1-p}+\lambda}{p}d^p(u_{i-1}^n,u^e)
\end{align*}
where in the first inequality we have again used minimality.
Let $\theta\in [0,1)$, divide by $1-\theta$, and take $\theta \to 1$ in order to get
\begin{align}
& E(\tau_i^n,u_{i-1}^n,u^e) +
  \frac{ (\tau_i^n)^{1-p}}{q}d^p(u_{i-1}^n,u^e)\nonumber\\
&\quad\leq  E(\tau_i^n,u_{i-1}^n,
  u_{i-1}^n)-\frac{\lambda}{p}d^p(u_{i-1}^n,u^e).\label{eq:qpo}
\end{align}
By taking the $q$-power of the slope estimate \eqref{eq:slopeestimate}
with $v=u^n_{i-1}$
we get 
$$ \sl^q(u^e) \leq (\tau_i^n)^{-p}d^p(u_{i-1}^n,u^e).$$
We use this to estimate from below the second
term on the left-hand side of \eqref{eq:qpo} obtaining 
\begin{align*}
& E(\tau_i^n,u_{i-1}^n,u^e) +\frac{\tau_i^n}{q}\sl^q(u^e)\leq  E(\tau_i^n,u_{i-1}^n,
  u_{i-1}^n) - \frac{\lambda}{p}d^p(u_{i-1}^n,u^e).% 
\end{align*}
As $E(\tau_i^n,u_{i-1}^n, u_{i-1}^n)=0$, given any $u_i^n\in M_G(\tau_i^n,u^n_{i-1})$ the latter
entails that  
\begin{align}
&G(\tau_i^n,u_{i-1}^n,u_i^n) \leq G(\tau_i^n,u_{i-1}^n,u^e) \nonumber\\
&\quad = E(\tau_i^n,u_{i-1}^n,u^e)
+\frac{\tau_i^n}{q}\sl^q(u^e) \leq -
  \frac{\lambda}{p}d^p(u_{i-1}^n,u^e).  \label{eq:han}
\end{align}
Recall now that the minimality $u^e\in M_E(\tau_i^n,u^n_{i-1})$ and
the nonnegativity of $\phi$
ensure  that
$$\frac{(\tau_i^n)^{1-p}}{p}d^p(u_{i-1}^n,u^e) \leq 
\phi(u_{i-1}^n).$$
Hence, inequality \eqref{eq:han} yields
\begin{equation}
  G(\tau_i^n,u_{i-1}^n,u_i^n)\leq 
\lambda^- (\tau_i^n)^{p-1}\phi(u_{i-1}^n).\label{eq:usero}
\end{equation}
Taking the sum on $i=1,\dots,m$ for $m\leq N^n$ we get
\begin{align*}
&\phi(u_m^n) + \frac1p \sum_{i=1}^m
(\tau_i^n)^{1-p}d^p(u_{i-1}^n,u_i^n) + \frac1q \sum_{i=1}^m\tau^n_i
\sl^q(u_i^n) -\phi(u^0) \\
&\quad = \sum_{i=1}^mG(\tau^n_i,u^n_{i-1},u^n_i) \leq \lambda^- (\tau^n)^{p-2}\sum_{i=0}^{m-1}
\tau^n_i\phi(u_i^n).
\end{align*}
We can hence use the discrete Gronwall Lemma and deduce that 
\begin{align*}
  &\phi(u_m^n) + \frac1p \sum_{i=1}^m
(\tau_i^n)^{1-p}d^p(u_{i-1}^n,u_i^n) + \frac1q \sum_{i=1}^m\tau^n_i
    \sl^q(u_i^n) \\
  &\quad \leq \phi(u^0) \,{\rm exp}\left( \lambda^- (\tau^n)^{p-2}
   t^n_m\right).
   \end{align*}
Going back to \eqref{eq:usero}, this entails that
$$(G(\tau_i^n,u_{i-1}^n,u_i^n))^+ \leq 
 \lambda^- (\tau_i^n)^{p-1}\phi(u^0)  \,{\rm exp}\left( \lambda^- (\tau^n)^{p-2}
  T\right). $$
Adding up for $i=1,\dots,N^n $ we get
\begin{align*}
  &
    \sum_{i=1}^{N^n} \big(G(\tau_i^n,u_{i-1}^n,u_i^n) \big)^+  \leq
\lambda^- (\tau^n)^{p-2}T\,\phi(u^0)  \,{\rm
    exp}\left( \lambda^- (\tau^n)^{p-2}T\right)=:R^n.
    \end{align*}
If $\lambda\geq 0$, we have that $R^n=0$ and condition
\eqref{eq:cond} trivially holds. If $\lambda<0$ and $p> 2$, one can
readily check that $R^n \to 0$ as $n\to \infty$ and 
\eqref{eq:cond} again holds.

%%%%%%%%%%%%%%%%%%%%%%%%%%%%%%%%%%%%%%%%%%%%%%%%%%%%%%%%%%%%%%%%%%%%%%%%%%%%%%%%%\
\section{Convergence without geodesic convexity}\label{sec:p3}
\setcounter{equation}{0}

We now turn to the proof of Theorem \ref{thm:main3}, where the
 convexity assumption is replaced by the generalized one-sided Taylor-expansion
assumption \eqref{eq:taylor}. The argument follows the general strategy
of \cite[Chap. 3]{Ambrosio08}, by revisiting the theory and adapting
it to the
incremental functional $G$ and to the case $p>1$. In particular, it is
fairly different with respect to that of Section \ref{sec:p2} and does
not rely on the existence of solutions of the Euler scheme. We prepare some
preliminary arguments in Subsections \ref{sec:meas}-\ref{sec:slope},
deduce an a priori estimate in Subsection \ref{sec:apriori} and eventually
present the proof of Theorem \ref{thm:main3}   in Subsection \ref{sec:conclude}.

\subsection{A measurable selection in $\tau \mapsto M_G(\tau,v)$}\label{sec:meas}
Let us recall that for all $\tau\in (0,\tau_*]$ and $v\in D(\phi)$ the
set of minimizers $M_G(\tau,v)$ is not empty. By additionally defining $M_G(0,v) =
\{v\}$, the set-valued function $\tau \in [0,\tau_*] \mapsto
M_G(\tau,v)$ has nonempty values. The aim of this section is to check
that it admits a measurable selection, namely,
\begin{equation}
  \label{eq:sel}
 \exists\, \tau \in [0,\tau_*] \mapsto u_\tau \in M_G(\tau,v) \ \ \text{measurable}.
\end{equation}

To this aim, we firstly check that $M_G(\tau,v)$ is closed for all
$\tau \in [0,\tau_*]$. Indeed, assume $\tau >0$ (the case $\tau=0$
being trivial) and let $u_k\in M_G(\tau,v)$ with $u_k \to 
u_\infty$. In particular, we have that 
$$\phi(u_k) + \frac{\tau^{1-p}}{p}d^p(v,u_k) +
\frac{\tau}{q}\sl^q(u_k) -\phi(v)= G(\tau,v,u_k) \leq G(\tau,v,w)$$
for any $w\in D(\sl).$  
Owing to the lower semicontinuity \eqref{eq:phicomp}-\eqref{eq:phisl}
we can pass to the lower limit and check that $G(\tau,v,u_\infty) \leq
G(\tau,v,w)$, so that $ u_\infty \in M_G(\tau,v)$ as well.  

Secondly, we check that $\tau \mapsto M_G(\tau,v)$ is measurable in
the sense of set-valued functions \cite{Wagner}. In particular, we
have to check that, for all $C \subset U$
closed, the set 
$$A=\{ \tau \in [0,\tau_*] \ : \ M_G(\tau,v)\cap C \not = \emptyset\}$$
is measurable. Indeed, one can prove that $A$ is closed: Take
$\tau_k \in A$ such that $\tau_k \to \tau_\infty$ and let $u_k \in
M_G(\tau_k,v)\cap C $. We have that 
\begin{align}
&\phi(u_k) + \frac{\tau^{1-p}_k}{p}d^p(v,u_k) +
\frac{\tau_k}{q}\sl^q(u_k) -\phi(v)= G(\tau_k,v,u_k) \nonumber\\
&\quad \leq G(\tau_k,v,v) = \frac{\tau_k}{q}\sl^q(v)<\infty.\label{eq:ttau}
\end{align}
One can hence deduce uniform estimates for $u_k$ and from compactness \eqref{eq:phicomp} one extracts a not
relabeled subsequence such that $u_k \to u_\infty$. If
$\tau_\infty>0$, by passing to the liminf in the minimality condition
for $u_k$ one gets
$$G(\tau_\infty,v,u_\infty)\leq \liminf_{k \to \infty}G(\tau_k,v,u_k)
\leq \liminf_{k \to \infty}G(\tau_k,v,w) = G(\tau_\infty,v,w) $$
for any $w\in D(\sl)$. This implies that $u_\infty\in
M_G(\tau_\infty,v)$. On the other
hand, if $\tau_\infty=0$ we obtain from \eqref{eq:ttau} that
$$ d^p(v,u_k)  \leq p\tau_k^{p-1}\phi(v)+\frac{p\tau_k^p}{q} \sl^q(v)\to 0,$$
 so that $u_\infty = v
\in M_G(0,v)$. Since $C$ is closed, $u_\infty\in C$ as well and
we have proved that $M_G(\tau_\infty,v) \cap C$ is not empty.
In particular,
$\tau_\infty\in A$ which is hence closed.

As the metric space $(U,d)$ is complete and separable and 
$\tau \mapsto M_G(\tau,v)$ has nonempty and closed values, the Ryll-Nardzewski
Theorem \cite{Ryll} applies and \eqref{eq:sel} holds. 

\subsection{Continuity of $\tau
\mapsto \haz G(\tau,v)$}\label{sec:diff}
We now turn our attention to the real map $\tau \in [0,\tau_*]
\mapsto \haz G(\tau,v)$ for some given $v\in D(\phi)$, where we
define $\haz G(0,v) = 0$. In order to check that this function is continuous on
$[0,\tau_*]$, take $\tau_k \in [0,\tau_*] \to \tau_\infty$ and $u_k
\in M_G(\tau_k,v)$. Following the argument of Subsection
\ref{sec:meas}, we can extract a not relabeled subsequence such that $u_k\to u_\infty \in M_G(\tau_\infty,v)$. 

If $\tau_\infty>0$
the lower semicontinuity \eqref{eq:phicomp}-\eqref{eq:phisl} implies
that
\begin{align*}
  &G(\tau_\infty,v,u_\infty) \leq \liminf_{k\to
  \infty} G(\tau_k,v,u_k) \leq \limsup_{k\to \infty} G(\tau_k,v,u_k)\\
&\quad 
  \leq \limsup_{k\to \infty} G(\tau_k,v,u_\infty) = G(\tau_\infty,v,u_\infty).
\end{align*}

The case
$\tau_\infty=0$ is even simpler as $u_\infty=v$ and we can compute
\begin{align}
 & 0=\haz G(0,v) = \phi(u_\infty) - \phi(v) \leq \liminf_{k \to \infty} \phi(u_k) - \phi(v) \leq
  \liminf_{k\to \infty} G(\tau_k,v,u_k) \nonumber\\
&\quad \leq \limsup_{k\to \infty}
  G(\tau_k,v,u_k)\leq \limsup_{k\to \infty} G(\tau_k,v,v) =
  \lim_{k\to \infty} \frac{\tau_k}{q} \sl^q(v) = 0. \label{eq:zero} 
\end{align}
In both cases, we have proved that $\haz G(\tau_k,v) \to \haz G(\tau_\infty,v)$.

\subsection{Differentiability of $\tau
\mapsto \haz G(\tau,v)$}\label{sec:diff} The aim of the subsection is
to show that $\tau
\mapsto \haz G(\tau,v)$ is even locally Lipschitz
continuous and to compute its almost-everywhere derivative, see equation
\eqref{eq:fallo3} below. 

Take $0<\tau_0 <\tau_1 <\tau_*$,
$u_0\in M_G(\tau_0,v)$, and $u_1\in M_G(\tau_1,v)$ where $v\in
D(\phi)$ is fixed. From minimality we
deduce
\begin{align*}
  \haz G(\tau_1,v) \leq G(\tau_1,v,u_0) = \haz G(\tau_0,v) +
  \frac{\tau_1^{1-p} - \tau_0^{1-p}}{p} d^p(v,u_0)
  +\frac{\tau_1 - \tau_0}{q}\sl^q(u_0)
\end{align*}
so that one has
$$ \haz G(\tau_1,v) - \haz G(\tau_0,v) \leq \frac{\tau_1^{1-p} - \tau_0^{1-p}}{p} d^p(v,u_0) +\frac{\tau_1 -
  \tau_0}{q}\sl^q(u_0).$$
By exchanging the roles of $\tau_0$ and $\tau_1$ we also get
$$ \haz G(\tau_0,v) - \haz G(\tau_1,v) \leq \frac{\tau_0^{1-p} - \tau_1^{1-p}}{p} d^p(v,u_1) +\frac{\tau_0 -
  \tau_1}{q}\sl^q(u_1).$$
By dividing by $\tau_1 - \tau_0$ we hence obtain 
\begin{align}
&\frac{\tau_1^{1-p} - \tau_0^{1-p}}{p(\tau_1 - \tau_0)}
                d^p(v,u_1)\leq \frac{\tau_1^{1-p} - \tau_0^{1-p}}{p(\tau_1 - \tau_0)}  d^p(v,u_1)
                +\frac{1}{q}\sl^q(u_1)\nonumber\\
  &\quad\leq \frac{\haz G(\tau_1,v) - \haz G(\tau_0,v)
  }{\tau_1- \tau_0}\nonumber\\
&\quad \leq \frac{\tau_1^{1-p} - \tau_0^{1-p}}{p(\tau_1 - \tau_0)}
                                 d^p(v,u_0) +\frac{1}{q}\sl^q(u_0)\leq
                                 \frac{1}{q}\sl^q(u_0).\label{eq:fallo}
\end{align}
The latter implies that $\tau
\mapsto \haz G(\tau,v)$ is locally
Lipschitz continuous on $(0,\tau_*]$. Indeed, take $0<\underline \tau< \tau_*$
and $\tau \in [\underline \tau,\tau_*]$.
Given $u_\tau \in M_G(\tau,v)$, we readily deduce   that 
\begin{align*}
  &d^p(v,u_\tau) \leq p  \tau^{p-1}_*\phi(v) + \frac{p
    \tau^p_*}{q}\sl^q(v), \\
  &\frac{1}{q}\sl^q(u_\tau) \leq
    \frac{1}{\underline \tau}\phi(v) + \frac{1}{q}\sl^q(v),\\
  &{}-\frac{\tau_1^{1-p} - \tau_0^{1-p}}{p(\tau_1 - \tau_0)}
   \leq \frac{1}{q\tau_0^p} \leq \frac{1}{q\underline \tau^p}.
    \end{align*}
In particular, moving from \eqref{eq:fallo}, for all $\underline \tau\in(0,\tau_*]$ we find $C_{\underline \tau}$
depending on $\underline \tau$, $\phi(v)$, and $\sl(v)$ such that 
$$\left| \frac{\haz G(\tau_1,v) - \haz G(\tau_0,v)
  }{\tau_1- \tau_0} \right|\leq C_{\underline \tau} \quad \forall \underline
  \tau< \tau_0 < \tau_1 < \tau_*.$$
Hence, $\tau \in (0,\tau_*]\mapsto \haz G(\tau,v)$ is locally
Lipschitz continuous and therefore almost
everywhere differentiable in $(0,\tau_*)$. 

Define now
\begin{align*}
  &\overline f(\tau_0,\tau_1) = \sup\left\{  \frac{\tau_1^{1-p} - \tau_0^{1-p}}{p(\tau_1 - \tau_0)}  d^p(v,u_1)
  +\frac{1}{q}\sl^q(u_1) \ : \ u_1\in M_G(\tau_1,v)   \right\},\\
&\underline f(\tau_0,\tau_1) = \inf\left\{  \frac{\tau_1^{1-p} - \tau_0^{1-p}}{p(\tau_1 - \tau_0)}  d^p(v,u_0)
  +\frac{1}{q}\sl^q(u_0) \ : \ u_0\in M_G(\tau_0,v)   \right\}.
\end{align*}
By using again relation \eqref{eq:fallo} one has
\begin{equation}
  \label{eq:fallo2}
  \ove f(\tau_0,\tau_1) \leq \frac{\haz G(\tau_1,v) - \haz G(\tau_0,v)
  }{\tau_1- \tau_0}\leq \underline f(\tau_0,\tau_1).
\end{equation}
Let $\tau \in (0,\tau_*)$ be such that $\tau \mapsto \haz G(\tau,v)$
is differentiable at $\tau$, take $h \in (0,\tau_*-\tau)$ and any
$u_{\tau+h} \in M_G(\tau+h,v)$. By arguing as in Subsection
\ref{sec:meas}, one can extract a not relabeled subsequence $u_{\tau+h}
\to u_\tau$ as $h \to 0$ and 
% From the minimality 
% $$ \phi(u_{\tau+h}) + \frac{(\tau+h)^{1-p}}{p}d^p(u,u_{\tau+h}) +
% \frac{\tau+h}{q}\sl^q(u_{\tau+h}) -\phi(u) \leq
% \frac{\tau+h}{q}\sl^q(u)$$
% we deduce that $d^p(u,u_{\tau+h})$, $\phi(u_{\tau+h})$, and
% $\sl^q(u_{\tau+h})$ are bounded, independently of $h$. By extracting a
% not relabeled subsequence and recalling that $\sigma$ is the topology
% induced by $d$ we obtain that 
% \begin{align*}
% &d^p(u,u_{\tau+h}) \to d^p(u,u_{\tau}), \\
% & \phi(u_\tau) \leq
% \liminf_{h\to 0} \phi(u_{\tau+h}), \\
% & \sl^q(u_\tau) \leq
% \liminf_{h \to 0} \sl^q(u_{\tau+h})
% \end{align*}
% for some $u_\tau$. Indeed, for all $w\in D(\sl)$ the above
% convergences  ensure that 
% $$ G(\tau,u,u_\tau) \leq \liminf_{h \to 0}G(\tau+h,u,u_{\tau+h}) \leq
% G(\tau+h,u,w) = G(\tau,u,w)$$
% so that 
check that
$u_\tau\in M_G(\tau,v)$. Moreover, going back to
\eqref{eq:fallo2} and choosing $\tau_0=\tau$ and $\tau_1=\tau+h$ we
deduce that 
\begin{align*}
 & -\frac{\tau^{-p}}{q}d^p(v,u_\tau)+\frac{1}{q}\sl^q(u_\tau)\leq \liminf_{h\to 0}\ove f(\tau,\tau+h)
  \\
&\quad= \frac{\rm d}{{\rm d} \tau} \haz G(\tau,v) \leq \liminf_{h\to 0}
  \underline  f(\tau,\tau+h) \leq
  -\frac{\tau^{-p}}{q}d^p(v,\tilde u_\tau)+\frac{1}{q}\sl^q(\tilde u_\tau)
\end{align*}
where $\tilde u_\tau$ is any element of $M_G(\tau,v)$. Passing to the
infimum in $M_G(\tau,v)$ left and right we get 
\begin{equation}\label{eq:fallo3}
\frac{\rm d}{{\rm d} \tau} \haz G(\tau,v) = \inf \left\{
    -\frac{\tau^{-p}}{q}d^p(v,u_\tau)+\frac{1}{q}\sl^q(u_\tau) \ : \
    u_\tau \in M_G(\tau,v) \right\}
\end{equation}
almost everywhere in $(0,\tau_*)$.

% Take now $\tau_0 \nearrow \tau_1$ left and $\tau_1\searrow \tau_0$
% right in order to conclude that 
%  \begin{align}
% &-\frac{1}{2\tau_1^2} d^2(u,u_1) +\frac{1}{2}\sl^2(u_1)\leq
%  \frac{\haz G(\tau_1,u) - \haz G(\tau_0,u) }{\tau_1- \tau_0} \nonumber\\
% &\quad \leq
%  -\frac{1}{2\tau_0^2} d^2(u,u_0) +\frac{1}{2}\sl^2(u_0). \label{eq:monotonie}
% \end{align}
% This proves that the map
% $$ \tau  \in (0,\tau_*) \mapsto -\frac{1}{2 \tau ^2} d^2(u,u_\tau )
% +\frac{1}{2}\sl^2(u_\tau ) \ \ \text{for} \ \ u_\tau  \in  M_G(\tau ,u)$$
% is nonincreasing. Moreover, by using that 
% $$ u_\tau \in M_G(\tau ,u) \ \Rightarrow \ \phi(u+\tau) + \frac{1}{2\tau}
% d^2(u,u_\tau) + \frac{\tau}{2}\sl^2(u_\tau) \leq
% \frac{\tau}{2}\sl^2(u),$$
% we get from \eqref{eq:monotonie} that 
% $$-\frac{1}{2\tau_1^2} d^2(u,u_1)  \leq
%  \frac{\haz G(\tau_1,u) - \haz G(\tau_0,u) }{\tau_1- \tau_0}.$$
% As both bounds above are independent of $\tau_0$, we deduce that  
%   $\tau \in (0,\tau_*) \mapsto \haz G(\tau,u)$ is
% locally Lipschitz and
% \begin{equation}
%   \label{eq:derivata}
% \frac{\rm d }{\rm d \tau} \haz G(\tau,u) =
% -\frac{1}{2\tau^2}d^2(u,u_\tau) + \frac12 \sl^2(u_\tau).
% \end{equation}

\subsection{Slope estimate}\label{sec:slope}
Let us prepare a version of the slope estimate
\eqref{eq:slopeestimate} adapted to our setting, namely for points in
$u \in M_G(\tau,v)$ for $v \in D(\phi)$ instead of $M_E(\tau,v)$. Let $w\in
D(\sl)$ be given. From minimality we
deduce that 
\begin{align*}
 & \phi(u) - \phi(w) + \frac{\tau}{q}\sl^q(u)  -
  \frac{\tau}{q}\sl^q(w) \leq
  \frac{\tau^{1-p}}{p}\left(d^p( v,w)-d^p(v,u)\right)\\
% & = \frac{1}{2\tau}\left(d( u^n_{i-1},w))-d(u^n_{i-1},u_\tau)\right)
% \left(d(w,u^n_{i-1}))+d(u^n_{i-1},u_\tau)\right)\\
&\leq d(u,w) \frac{\tau^{1-p}}{p}\sum_{k=1}^\infty \binom{p}{k}d^{k-1}( u,w)d^{p-k}(v,u).
\end{align*}
By assuming that $w\not = u$, dividing by $d(u,w) $, and taking
$w \to u$ we get 
\begin{equation}
  \label{eq:slopeestimate2}
  |\partial ( \phi + \tau\sl^q/q)|(u) \leq \tau^{1-p}
  d^{p-1}(v,u)\quad \forall u \in M_G(\tau,v).
\end{equation}
This proves in particular the additional regularity 
$$ M_G(\tau,v) \subset D(\partial ( \phi + \tau\sl^q/q))$$
for minimizers of $G$.

\subsection{A priori estimate}\label{sec:apriori}
Let now $\{u_i^n\}_{i=0}^{N^n}$ solve the incremental minimization problem
\eqref{eq:min} with $u^0$ replaced by the approximating $u^{0n}\in M_E(\tau^n,u^0)$. From minimality we obtain that
\begin{equation} \phi(u^n_i) +
\frac{\tau^n_i}{q}\sl^q(u^n_i) + \frac1p (\tau^n_i)^{1-p} d^p(u^n_{i-1},u^n_i) \leq \phi(u^n_{i-1})   +
\frac{\tau^n_i}{q}\sl^q(u^n_{i-1}).\label{eq:3}
\end{equation}
Taking into account the one-sided nondegeneracy of the time partition $$(\tau^n_i
-\tau^n_{i-1})^+/\tau^n_{i-1} \leq  C \tau^n$$ we can control the above
right-hand side of \eqref{eq:3}  as follows
\begin{align*}
  &\phi(u^n_{i-1})   +
\frac{\tau^n_i}{q}\sl^q(u^n_{i-1}) =
 \phi(u^n_{i-1})   +
 \frac{\tau^n_{i-1}}{q}\sl^q(u^n_{i-1}) + \frac{\tau_i^n - \tau^n_{i-1}}{q}\sl^q(u^n_{i-1}) \\
 &\quad \leq \phi(u^n_{i-1})   +
\frac{\tau^n_{i-1}}{q}\sl^q(u^n_{i-1}) +  C\tau^n
   \frac{\tau^n_{i-1}}{q}\sl^q(u^n_{i-1}) .
   \end{align*}
Owing to this bound, we can take the sum in \eqref{eq:3} for $i=2,\dots,m$ and get
\begin{align*}
&\phi(u_m^n) + \frac{\tau^n_m}{q}\sl^q(u^n_m) +
  \frac1p\sum_{i=1}^{m}(\tau^n_i)^{1-p}d^p(u_{i-1}^n, u_i^n) \nonumber\\
  &\quad \leq \phi(u^n_{1})   +
\frac{\tau^n_{1}}{q}\sl^q(u^n_{1}) +  C\sum_{i=2}^{m}\tau^n
    \frac{\tau^n_{i-1}}{q}\sl^q(u^n_{i-1})
     \nonumber\\
  &\quad \stackrel{\eqref{eq:3}}{\leq} \phi(u^{0n})   +
\frac{\tau^n }{q}\sl^q(u^{0n}) +  C\sum_{j=1}^{m-1}\tau^n
    \frac{\tau^n_{j}}{q}\sl^q(u^n_{j}).
\end{align*}
By applying the discrete Gronwall Lemma we hence obtain  
\begin{align}
  &\phi(u_m^n) + \frac{\tau^n_m}{q}\sl^q(u^n_m) +
  \frac1p\sum_{i=1}^{m}(\tau^n_{i})^{1-p}d^p(u_{i-1}^n, u_i^n)
  \nonumber\\
  &\quad \leq C\left(\phi(u^{0n})   +
\frac{\tau^n }{q}\sl^q(u^{0n}) \right).
  \label{eq:4}
\end{align}
Recall now that $u^{0n} \in M_E(\tau^n,u^0)$ and use the slope
estimate \eqref{eq:slopeestimate} to get that
$$\frac{\tau^n }{q}\sl^q(u^{0n}) \leq \frac1q (\tau^n)^{1-p}
d^p(u^0,u^{0n}) \leq \frac{p}{q}\phi(u^0).$$
Hence, $\{u^{0n}\}$ are in particular $d$-bounded and the bound \eqref{eq:4} entails the estimate
\begin{align}
  &\phi(u_m^n) + \frac{\tau^n_m}{q}\sl^q(u^n_m) +
  \frac1p\sum_{i=1}^{m}\tau^{1-p}d^p(u_{i-1}^n, u_i^n) \leq
  C\left(1+\frac{p}{q}\right) 
    \phi(u^0) \nonumber\\
  &\quad \forall m =1, \dots, N^n, \ \forall n. \label{eq:estimate}
  \end{align}

\subsection{Conclusion of the proof}\label{sec:conclude}
For all $i =1,\dots,N^n$ and $\tau_0 \in (0,\tau^n_i]$
we use the Lipschitz continuity of $\tau \in  (0,\tau^n_i]\mapsto \haz
G(\tau,u^n_{i-1})$ and write   
\begin{align}
  \haz G(\tau_i^n,u_{i-1}^n) =   \haz
  G(\tau_0,u_{i-1}^n) +\int_{\tau_0}^{\tau_i^n}\frac{\rm d
  }{\rm d \tau} \haz G(\tau,u^n_{i-1})\, {\rm d} \tau . \label{eq:eatme}
\end{align}

Let now $\tau\in [\tau_0,\tau^n_i] \mapsto u_\tau $ be a measurable selection in
$M_G(\tau,u^n_{i-1})$. The existence of such a selection is ascertained in Subsection \ref{sec:meas}.
Take $\tau_0 \to 0$ in \eqref{eq:eatme}, use $\haz
G(\tau_0,u^n_{i-1}) \to 0$ from 
\eqref{eq:zero} and \eqref{eq:fallo3} to get 
\begin{equation}
G(\tau^n_i,u^n_{i-1},u^n_i) \leq  \int_0^{\tau^n_i}
\left( -  \frac{\tau^{-p}}{q}
  d^p(u^n_{i-1},u_\tau)+\frac1q\sl^q(u_\tau)\right){\rm d} \tau.\label{eq:controlG}
\end{equation}
%for $\tau \mapsto u_\tau \in M_G(\tau, u^n_{i-1})$. % \UUU qui manca la misurabilita`\EEE

In order to conclude the proof of Theorem \ref{thm:main2}, one has to check that condition
\eqref{eq:cond} holds, so that Theorem \ref{thm:main1} applies. This calls
for controlling the right-hand side of \eqref{eq:controlG}. By means
of the slope estimate \eqref{eq:slopeestimate2}  for $v=u^n_{i-1}$
and $u=u_\tau$ we can control the right-hand of \eqref{eq:controlG} as
$$
G(\tau^n_i,u^n_{i-1},u^n_i) \leq   \int_0^{\tau^n_i}
\left(\frac1q\sl^q(u_\tau) -  \frac1q  |\partial ( \phi +
  \tau\sl^q/q)|^q(u_\tau) \right){\rm d} \tau.
$$
We now use estimate \eqref{eq:estimate} and the generalized one-sided Taylor expansion condition
\eqref{eq:taylor} in order to conclude that 
\begin{align*}
  &\sum_{i=1}^{N^n} \left( G(\tau^n_i,u^n_{i-1},u^n_i)\right)^+ \leq \frac1q
  \sum_{i=1}^{N^n} \int_0^{\tau^n_i}g(\tau)\,{\rm d} \tau \\
&\quad = \frac1q
  \sum_{i=1}^{N^n} \tau^n_i  
  \left(\frac{1}{\tau^n_i}\int_0^{\tau^n_i}g(\tau)\,{\rm d} \tau
  \right) \leq \frac{T}{q}  \frac{1}{\tau^n}\int_0^{\tau^n}g(\tau)\,{\rm d} \tau. 
\end{align*}
As $(1/\tau^n)\int_0^r g(\tau^n)\, {\rm d} \tau \searrow 0$ as
$\tau^n\to 0$ condition \eqref{eq:cond} holds. The statement hence
follows from Theorem \ref{thm:main1}.

%%%%%%%%%%%%%%%%%%%%%%%%%%%%%%%%%%%%%%%%%%%%%%%%%%%%%%%%%%%%%%%%%%%%%%%%%%%%%%%%%\
\section{Applications in linear spaces}\label{sec:appl0}

We collect in this section some comments on the application of the
abstract convergence results of Theorem
\ref{thm:main1}-\ref{thm:main3} in linear finite and
infinite-dimensional spaces.

Let us start from the convex case of Theorem \ref{thm:main2}. We hence
restrict to $p=2$,  for assumption \eqref{eq:F2} cannot hold for $p\not
=2$ in linear spaces, as commented in Subsection
\ref{sec:assumptions}. Correspondingly, the potential  $\phi$ is requires to be
convex ($\lambda \geq 0$). 

In the finite-dimensional ODE case, let the proper, convex potential
$\phi:\Rz^d \to  [0,\infty]$ and the initial datum $u^0 \in
D(\phi)$ be given. In
this case, we have that $\sl(u) = |(\partial \phi (u))^\circ|$, where
$(\partial \phi (u))^\circ$ is the element of minimal norm in the
convex and closed set $\partial \phi (u)$. In particular, $\sl(u)$ is lower
semicontinuous. As such, the new minimizing-movements scheme
\eqref{eq:min} has a solution $\{u_i^n\}$ for any partition and the
corresponding interpolants
 converge to a solution
of  $u' +\partial \phi(u) \ni 0$, up to subsequences.

In order to give an application of Theorem \ref{thm:main2} in infinite dimensions, we consider \begin{equation}
  \label{eq:pde}
   \partial_t u - \nabla {\cdot} \beta (\nabla
  u) + \alpha (u)\ni 0 \quad \text{in}\  \ \Omega \times (0,T).
\end{equation}
Here, $\Omega \subset \Rz^d$ is open, bounded, and smooth, $u: \Omega \times
(0,T) \to \Rz$ is scalar-valued, and $\partial_t$ and $\nabla$ indicate partial
derivatives in time and space, respectively. We assume that $\beta =
\partial \haz \beta$ and $\alpha = \partial  \alpha$ where the potentials $\haz
\beta: \Rz^d\to [0,\infty]$ and $\haz \alpha:\Rz \to [0,\infty]$ are
proper and convex. In addition, we assume $\haz \beta$ to be coercive
in the following sense 
\begin{equation} \label{beta}
\exists \, c_\beta >0, \ m>\frac{2d}{d+2} : \quad \haz \beta (\xi) \geq c_\beta |\xi|^m -
\frac{1}{c_\beta} \quad \forall \xi \in \Rz^d.
\end{equation}
Equation \eqref{eq:pde} is intended to be complemented with homogeneous
Dirichlet boundary conditions (other choices being
of course possible) hence corresponding to the gradient flow in
$U=L^2(\Omega)$ of the functional 
\begin{equation}
\phi(u)= 
\left\{
  \begin{array}{ll}
\disp\int_\Omega \big(\haz \beta (\nabla u )  + \haz \alpha (u) \big) \,
    \dx  &\text{for} \ \ u \in W^{1,m}_0(\Omega), \  \\
&\quad\text{with} \  \haz \beta (\nabla u )  + \haz \alpha
  (u)\in L^1(\Omega)   \\[2mm]
\infty &\text{elsewhere in}   \ \ L^2(\Omega). 
  \end{array}
\right.\label{def:phi}
   \end{equation}
As $\phi:U \to [0,\infty]$ is convex, proper, and lower
semicontinuous, we have that $u \mapsto \partial \phi(u)$ is
strongly-weakly closed and  
 $\sl(u) = \| (\partial \phi(u))^\circ\|$
(norm in $L^2(\Omega)$) is lower semicontinuous. Moreover, $\partial
\phi$ fulfills the chain rule \cite[Lem. 3.3]{Brezis73}, so that $\sl$ is a
strong upper gradient. Note that the sublevels of $\phi$ are bounded
in $W^{1,m}_0 (\Omega)$, which embeds compactly into $L^2(\Omega)$. We
can hence apply Theorem \ref{thm:main2}. In particular, the new minimizing-movements scheme
\eqref{eq:min} has a solution, which  
 converges to a solution
of \eqref{eq:pde}, up to subsequences.

Let us now turn to some application of Theorem \ref{thm:main3} to
nonconvex problems. In the finite-dimensional case,  assume $\phi$ to
be twice differentiable and coercive with $\nabla \phi$ and ${\rm
  D}^2\phi$ locally bounded. Then, one computes
\begin{align*}
  &|\nabla \phi(u)|^q - |\nabla (\phi(u) + \tau |\nabla
    \phi(u)|^q/q)|^q \nonumber\\
    &\quad =
   |\nabla \phi(u)|^q - |\nabla \phi(u) + \tau |\nabla \phi(u)|^{q-2}
  {\rm D}^2\phi(u) \, \nabla \phi(u)|^q \\
  &\quad \leq  \left(\big|\nabla \phi(u)+\tau |\nabla \phi(u)|^{q-2}
  {\rm D}^2\phi(u) \, \nabla \phi(u)\big| + \big| \tau |\nabla \phi(u)|^{q-2}
    {\rm D}^2\phi(u) \, \nabla \phi(u)\big|\right)^q \nonumber\\
  &\qquad - \big|\nabla \phi(u) + \tau |\nabla \phi(u)|^{q-2}
    {\rm D}^2\phi(u) \, \nabla \phi(u)\big|^q \nonumber\\
  &\quad \leq \tau \sum_{k=1}^\infty\binom{q}{k}\big|\nabla \phi(u)+\tau |\nabla \phi(u)|^{q-2}
  {\rm D}^2\phi(u) \, \nabla \phi(u)\big|^{q-k} \big|\nabla \phi(u)|^{q-2}
    {\rm D}^2\phi(u) \, \nabla \phi(u)\big|^k\nonumber\\
  &\quad \leq \tau \left( \big|\nabla \phi(u)+\tau |\nabla \phi(u)|^{q-2}
  {\rm D}^2\phi(u) \, \nabla \phi(u)\big| + \big||\nabla \phi(u)|^{q-2}
    {\rm D}^2\phi(u) \, \nabla \phi(u)\big|\right)^q.
\end{align*}
Hence, the one-sided Taylor-expansion condition \eqref{eq:taylor} holds for the choice
$$g(\tau) = \tau \sup_{\phi(v)\leq C}\left( |\nabla \phi(v)|
  + 2  |\nabla \phi(v)|^{q-1} |{\rm D}^2 \phi(v)|\right)^q.$$
Note that the above computation simplifies in case $p=2$, for
we have
\begin{align*}
 & |\nabla \phi(u)|^2 - |\nabla (\phi(u) + \tau |\nabla
    \phi(u)|^2/2)|^2  = |\nabla \phi(u)|^2 - |\nabla \phi(u) + \tau  
  {\rm D}^2\phi(u) \, \nabla \phi(u)|^2 \nonumber\\
  &\quad - \tau^2 |{\rm D}^2\phi(u) \, \nabla \phi(u)|^2 - 2\tau
    \nabla \phi(u) {\cdot} ({\rm D}^2\phi(u) \, \nabla \phi(u)).
\end{align*}
In particular, if $\phi$ is convex condition \eqref{eq:taylor} holds
with the trivial choice $g(\tau)=0$. In all cases, if ${\rm D}^2\phi(u) $  is bounded below on sublevels of $\phi$ in the following
sense
\begin{equation}
  \forall C>0, \, \exists c>0, \, \forall v,\, \xi \in \Rz^d \
  \text{with} \ \phi(v)\leq C: \quad \xi {\cdot} {\rm
  D}^2\phi(v)\xi  \geq -  \frac{c}{2} | \xi|^2\label{anche}
\end{equation}
and $\nabla \phi(u)$ is
bounded on the sublevels of $\phi$, namely, $|\nabla \phi(u)| \leq \ell
(\phi(u))$ for some $\ell$ increasing,  we can choose $g(\tau)= 2\tau
c (\ell(C))^2$ 
in order to get again condition \eqref{eq:taylor}. This in particular
applies to $\phi\in C^2$ and coercive. In all cases, we can apply
Theorem \ref{thm:main3} and deduce that the solution of the new
minimizing-movements scheme converges up to subsequences to a solution
of \eqref{ode}.

Let us now turn to the infinite-dimensional case. To simplify
notation, let again $p=2$ and $U=L^2(\Omega)$ (the case $p\not =2$ and
$U=L^p(\Omega) $ can
also be treated) and define $\phi$ as in \eqref{def:phi} by
dropping the convexity requirement on   $\haz
\alpha$. More precisely, we ask $\beta = {\rm D} \haz \beta \in C^2(\Rz^d;\Rz^d)$ and
$\alpha = \haz \alpha' \in C^2(\Rz)$ and $\haz \beta$ fulfill the
coercivity \eqref{beta}.
% $$\exists \, C>0, \  m>\frac{2d}{d+2}: \quad \haz \beta(\xi) \geq \frac{1}{C}|\xi|^m
% - C \quad \text{and} \quad |\beta(\xi)| \leq c(1+
% |\xi|^{m-1}) \quad \forall \xi \in \Rz^d.$$
% In order to further simplify notation, assume additonally that $\haz \alpha(u) \leq
% C(1+|u|^{m_*})$ for all $u\in \Rz$, where $m_*=md/(d-m)^+$. 
In this
case, we have that 
\begin{align*}&\partial \phi(u)=  - \nabla {\cdot} \beta(\nabla u) +
  \alpha(u),\\
  &\quad 
\text{with} \ \ D(\partial \phi)= \{u \in L^2(\Omega) \ : \  - \nabla{\cdot}  \beta(\nabla u) + \alpha(u) \in L^2(\Omega) \}.
\end{align*}
Recall that the {\it Fr\'echet subdifferential} \cite{Rockafellar}
of $\psi:U \to [0,\infty]$ at $u\in D(\psi)$ is the set  
$$
\partial \psi(u)=\left\{ \xi \in U \ : \ \liminf_{v \to u}
  \frac{\psi(v)- \psi(u)- (\xi,v-u)}{\|v-u\|} \geq 0\right\}
$$
and $D(\partial \psi) = \{u \in D(\psi) \ : \ \partial \psi(u)\not = \emptyset\}$.
In case of   $\psi(u)=
\| \partial \phi(u)\|^2/2$ we obtain that the Fr\'echet
subdifferential is single-valued and
 \begin{align*}&   \partial \frac12 \| \partial \phi(u)\|^2  = \nabla {\cdot} \left( \DD \beta(\nabla
    u) \nabla \left( \nabla {\cdot} \beta(\nabla u) -
    \alpha(u)\right)\right) - \left(\nabla {\cdot} \beta(\nabla u) -
                 \alpha(u) \right) \alpha'(u)
 \end{align*}
with domain given by   
 \begin{align}    D\left( \partial \frac12 \| \partial
                                       \phi(u)\|^2 \right) &= \Big\{u \in
                                       D(\partial \phi)\ : \
                                       \partial   \| \partial
                 \phi(u)\|^2 \in L^2(\Omega), \nonumber\\
   &\quad \text{with} \ \ \left( \nabla {\cdot} \beta(\nabla u) -
    \alpha(u)\right)  \DD \beta(\nabla
    u)  \nu=0 \ \ \text{on} \ \ \partial \Omega  \Big\}.\label{eq:domain}
 \end{align}
 In particular, an extra natural boundary condition arises, where $\nu$
 denotes the outer normal vector to $\partial \Omega$.
 In the linear case of $\beta(\xi)=\xi$ (see
 Subsection \ref{sec:illu}), we have that $\DD \beta = I$ (identity
 matrix) and we deduce again
 \begin{align*}
  & \partial \phi(u) = -\Delta u, \quad D(\partial \phi) = \{ u \in
   H^1_0(\Omega) \ : \ -\Delta u \in L^2(\Omega)\}=H^2(\Omega) \cap
   H^1_0(\Omega), \\
   &\partial \frac12\| \Delta u\|^2 = \Delta^2 u, \\
   &D\left( \partial \frac12\| \Delta u\|^2\right) = \{ u \in
  H^2(\Omega) \cap
   H^1_0(\Omega) \ : \  \Delta^2 u \in L^2(\Omega) \ \text{and} \
     \Delta u =0 \ \text{on} \ \partial \Omega\}\\
   &\quad =\{  u \in
  H^4(\Omega) \cap
   H^1_0(\Omega) \ : \  \Delta u \in H^2(\Omega) \cap
   H^1_0(\Omega)\}.
 \end{align*}

 In order to assess the one-sided Taylor-expansion condition
 \ref{eq:taylor} we argue as follows
 \begin{align}
   & \sl^2(u) - |\partial (\phi + \tau \partial \sl^2/2)|^2(u) \nonumber\\
   & \quad=\|\partial \phi(u)\|^2 - \| \partial \phi(u) + \tau \partial \|
     \partial \phi(u)\|^2/2\|^2 \nonumber\\
   &\quad = - \tau^2 \|  \nabla {\cdot} \left( \DD \beta(\nabla
    u) \nabla \left( \nabla {\cdot} \beta(\nabla u) -
    \alpha(u)\right)\right) - \left(\nabla {\cdot} \beta(\nabla u) -
     \alpha(u) \right) \alpha'(u) \|^2 \nonumber\\
   &\qquad + 2\tau \int_\Omega \Bigg(\nabla {\cdot} \left( \DD \beta(\nabla
    u) \nabla \left( \nabla {\cdot} \beta(\nabla u) -
    \alpha(u)\right)\right) - \left(\nabla {\cdot} \beta(\nabla u) -
     \alpha(u) \right) \alpha'(u) \Bigg){\cdot} \nonumber\\
   &\qquad \qquad \qquad \qquad {\cdot} \left( \nabla {\cdot} \beta(\nabla u) -
     \alpha(u)\right)\, \dx \nonumber\\
   &\quad \leq -2\tau \int_\Omega  \nabla\left( \nabla {\cdot} \beta(\nabla u) -
    \alpha(u)\right){\cdot}\DD \beta(\nabla
    u) \nabla \left( \nabla {\cdot} \beta(\nabla u) -
     \alpha(u)\right)\, \dx \nonumber\\
   &\qquad -2\tau \int_\Omega \alpha'(u)  \left( \nabla {\cdot} \beta(\nabla u) -
     \alpha(u)\right)^2\, \dx\label{eq:local}
 \end{align}
 where we have used also the additional natural condition from
 \eqref{eq:domain} in the last inequality. The one-sided
 Taylor-expansion condition \eqref{eq:taylor} then holds if $\haz \beta$ and $\haz
 \alpha$ are convex. 

In addition, some nonconvex   $\haz \alpha$ can be
 considered as well. 
 Assume $m>d$. %, let $\DD\beta $ and $\haz \alpha$  be unformly bounded from below in the following sense
 % \begin{align*}
 %  &\exists c_\alpha>0: \quad  \haz \alpha(u) \geq -c_\alpha (1+|u|^2) \ \forall u \in \Rz.
 % \end{align*} 
Due to the coercivity of $\beta$, one has that the sublevels of $\phi$
 are bounded in $W^{1,m}$ hence in 
 $L^\infty$. In particular, $\phi(u) \leq c \ \Rightarrow \ \|u \|_{L^\infty} \leq\ell (c)$
 for some $\ell: (0,\infty)\to (0,\infty)$ increasing. Assume 
 $u^0$ to be given and use \eqref{eq:estimate} to bound
 $\phi(u)$. Owing to the above discussion we hence have that $\| u
 \|_{L^\infty} \leq \ell (2C\phi(u^0))$ along the discrete evolution, where $C$ is the constant in
 \eqref{eq:estimate}. 
Let now $C_{\rm P}>0$ be the Poincar\'e constant
 giving $\| w \|_{L^2}^2 \leq C_{\rm P}\| \nabla w\|_{L^2}^2$
 for all $w \in H^1_0(\Omega)$. Assume $\haz \alpha$ to be such that
 $\alpha'$ locally bounded from below. Under the following smallness assumption 
 $$\inf\big\{\alpha'(r) \ : \ |r| \leq \ell (2C\phi(u^0)) \big\}\geq -
 \frac{c_\beta}{C_{\rm P}}$$
 one has that the right hand side of \eqref{eq:local} can be
 controlled from above as
 \begin{align*}
  & -2\tau c_\beta \| \nabla \left( \nabla {\cdot} \beta(\nabla u) -
     \alpha(u)\right)\|^2_{L^2}+2\tau c_\beta |\Omega|+2\tau \frac{c_\beta}{C_{\rm P}} \|\nabla {\cdot} \beta(\nabla u) -
     \alpha(u)\|^2_{L^2}
 \end{align*}
 and the one-sided Taylor-expansion condition
 \eqref{eq:taylor} follows with $g(\tau)=2\tau c_\beta |\Omega|$, at
 least on the relevant energy sublevel. In this case, Theorem
 \ref{thm:main3} again ensures that the solution of the new
 minimizing-movement scheme converges to a solution of \eqref{eq:pde}, up
 to subsequences.

\section{Applications in Wasserstein spaces}\label{sec:Wass}

Let us now give some detail in the direction of the application of the
above theory to the case of the nonlinear diffusion equation
\eqref{eq:introd2}. To start with, let us specify the space of
probability measures of finite $p$-moment as 
$$U = \mathcal P_p(\Rz^d) = \left\{ u \in \mathcal P(\Rz^d) \ : \
  \int_{\Rz^d}|x|^p{\rm d}u(x)<+\infty \right\}$$
where $\mathcal P(\Rz^d)$ denotes probability measures on $\Rz^d$, and
endow it with the $p$-Wasserstein distance
$$W_p^p(u_1,u_2) = \inf\left\{ \int_{\Rz^d \times \Rz^d}|x-y|^p{\rm d}
  \mu (x,y) \ : \ \mu \in \mathcal P(\Rz^d {\times} \Rz^d), \ \pi^1_\#
  \mu = u_1, \pi^2_\#
  \mu = u_2\right\}$$
where   $u_1, \, u_2 \in \mathcal P_p(\Rz^d)$ and $\pi^i_\# $ denotes
the push-forward of the projection $\pi^i$ on the $i$-th
component. Let $\sigma$ indicate
the {\it narrow} topology, namely, $u_n \stackrel{\sigma}{\to} u$ iff
$$\lim_{n \to \infty}\int_{\Rz^d} f(x) \, {\rm d} u_n(x) \to
\int_{\Rz^d} f(x) \, {\rm d} u(x)\quad \forall f: \Rz^d \to \Rz \
\text{continuous and bounded}.$$
Note that $(\mathcal P_p(\Rz^d) , W_p)$ is a complete
metric space \cite[Prop.~7.1.5]{Ambrosio08} and that $\sigma$ is
compatible with $W_p$ \cite[Lemma~7.1.4]{Ambrosio08}, namely,
assumptions \eqref{eq:X}-\eqref{eq:compat} hold.

Let su now fix some assumptions on 
potentials $V$, $F$, and $W$. We follow the setting of
\cite[Sec.~10.4.7]{Ambrosio08}, also referring to \cite[Sec.~7]{rsss}
for some additional discussion. In particular, we assume 
\begin{align}
  &V:\Rz^d \to [0,\infty) \ \ \text{$(\lambda,2)$-convex with} \ \limsup_{|x| \to
  \infty}\frac{V(x)}{|x|^2}=\infty, \label{eq:V}\\[1mm]
  &F:[0,\infty) \to \Rz \ \ \text{convex, differentiable, superlinear for $|x|\to \infty$, 
    $F(0)=0$, and} \ \nonumber\\
  &\quad \exists C_F>0: \quad F(x+y) \leq C_F (1+F(x)+F(y)) \quad \forall x,\, y \in
    \Rz^d,\nonumber\\
 &\quad r \in (0,\infty)\mapsto r^d F(r^{-d}) \ \ \text{is convex and nonincreasing}, \label{eq:F}\\[3mm]
 & W: \Rz^d  \to  [0,\infty) \ \ \text{convex, differentiable, even,
   such that} \ \nonumber\\
  &\quad \exists C_W>0:\quad W(x+y) \leq C_W (1+W(x)+W(y)) \quad \forall x,\, y \in
    \Rz^d.\label{eq:W}
\end{align}
Note that the assumptions on $F$ cover the classical cases $F(r) = r
\ln r$ and $F(r) = r^m$ for $m>1$, respectively related to
Fokker-Planck and porous media equations.

Under assumptions
\eqref{eq:V}-\eqref{eq:W} we have that the potential $\phi$ from
\eqref{eq:f} is $(\lambda,2)$-geodesically convex. Combining this with
the $(1,2)$-generalized-geodesic convexity of $u \mapsto
W_2^2(v,u)$ \cite[Lemma 9.2.1]{Ambrosio08} one has that condition
\eqref{eq:F2} holds. Note that resorting to generalized-geodesic convexity is here
crucial,  for the  
Wasserstein space $({\mathcal P}_2(\Rz^d), W_2)$ is 
positively curved \cite[Prop.~3.1]{Ohta}, namely, $u \mapsto
W_2^2(v,u)$ is actually $(1,2)$-geodesically {\it concave}. 
In addition, $\phi$
has $\sigma$-sequentially compact
sublevels and its local slope $\sl$ is a
strong upper gradient and is $\sigma$-sequentially lower
semicontinuous  \cite[Prop.~10.4.14]{Ambrosio08}. In particular,
\eqref{eq:phicomp}-\eqref{eq:phisl} holds and we have the
following.

\begin{proposition}\label{prop:prop2}
  Assume \eqref{eq:V}-\eqref{eq:W} and $u^0\in \mathcal
P_2(\Rz^d)$ with $\phi(u^0)<\infty$. Let  
  $\{0=t_0^n<t_1^n<\dots<t^{n}_{N^n}=T\}$ be a sequence of partitions
  with $\tau^n:=\max (t^n_i-t^n_{i-1}) \to 0$ as $n\to
  \infty$. Moreover, let $u^n_i \in M_G(\tau^n_i,u^n_{i-1})$ for $i=1,
  \dots, N^n$. Then, up to a not relabeled subsequence, we have that
  $\ove u^n(t) \sto u(t)$, where $u \in AC^2([0,T];\mathcal
  P_2(\Rz^d))$ and there exists a density $\rho: t \in [0,T] \to
  L^1(\Rz^d)$ such that 
$u(t) = \rho(t) \mathcal L^d$, $
\int_{\Rz^d}\rho(x,t)\, {\rm d} \mathcal L^d(x)=1$, and $
\int_{\Rz^d}|x|^2\rho(x,t)\, {\rm d} \mathcal L^d(x)<\infty$  for
  all $t \in [0,T]$, satisfying $u^0=\rho(\cdot,0)\mathcal L^d$ and
the nonlinear diffusion equation
$$ \partial_t \rho - {\rm div}\big(   \rho\nabla ( V +
    F'(\rho)
 +      W \ast \rho)   \big)=0\quad \text{in} \ \ \mathcal D'(\Rz^d
\times (0,T)).$$
\end{proposition}

Let us now turn to an application of the one-sided Taylor-expansion
condition \eqref{eq:taylor} for general $p$. In the metric situation of \eqref{eq:f},
one can use such
condition  in the purely trasport case $F=0$ and $W=0$. By assuming periodic boundary conditions,
we formulate the problem on the torus
$\Tz^d=\Rz^d/\Zz^d$. Let $V \in C^3(\Tz^d)$ and define
\begin{equation}
  \label{eq:ff}
  \phi(u) = \int_{\Tz^d} V(x)\, {\rm d} u(x) \quad \forall u \in
  \mathcal P(\Tz^d).
\end{equation}
From \cite[Prop.~10.4.2]{Ambrosio08} we have 
that
$$\sl^q(u) = \int_{\Tz^d} |\nabla V(x)|^q  {\rm d} u(x).$$
In case $p=2$ we obtain  
$$ \phi(u) + \frac{\tau}{2}\sl^2(u)
=\int_{\Tz^d}\left(V(x)+\frac{\tau}{2}|\nabla V(x)|^2 \right) 
{\rm d} u(x) =: \int_{\Tz^d} \tilde V(x)  \,
{\rm d} u(x).$$
One readily checks that ${\rm D}^2\tilde V = {\rm D}^2 V +\tau {\rm
  D}^3 V \,\nabla V + \tau  {\rm D}^2 V  \,{\rm D}^2 V $ is bounded below. We can
hence  apply \cite[Prop.~10.4.2]{Ambrosio08} once more and deduce
that
\begin{align*}
 &\quad  \sl^2(u) - |\partial(\phi + \tau \sl^2/2)|^2(u) = \int_{\Tz^d}
  \left(|\nabla V(x)|^2 - |\nabla \tilde V(x)|^2 \right) {\rm d}
  u(x) \\
  &\quad =\int_{\Tz^d}
  \left(|\nabla V(x)|^2 - |\nabla  V(x) + \tau {\rm D}^2V(x) \nabla V(x)|^2\right) {\rm d}
    u(x)\\
  &\quad = - \int_{\Tz^d} \tau^2| {\rm D}^2V(x) \nabla V(x)|^2\,  {\rm d}
  u(x) -  2\int_{\Tz^d} \tau \nabla V(x) {\cdot} {\rm D}^2V(x) \nabla
    V(x)  \,{\rm d}
    u(x) \\
  &\quad \leq 2\tau \lambda^- \int_{\Tz^d}|\nabla V(x)|^2\, {\rm d}u(x) \leq
    2\tau \lambda^- \| \nabla V\|_{L^\infty(\Tz^d)}^2
\end{align*}
where we have defined $\lambda = \min\{ \xi {\cdot} {\rm
  D}^2V(x) \xi \ : \ x \in \Tz^d,\, \xi \in \Rz^d, \, |\xi|=1\}$. Hence, condition \eqref{eq:taylor} holds  with $g(\tau):=2\tau \lambda^-
\| \nabla V\|_{L^\infty(\Tz^d)}^2$ (and, in particular, $g(\tau)=0$ if $V$ is convex).

In fact, the above computation
can be adapted to the case $p\not = 2$ by letting $\tilde V = V +\tau
|\nabla V|^q/q$. Let us shorten notation by denoting by
$\xi(x) = \nabla V (x)$ and by $A(x) = {\rm D}^2 V(x)$. Then, 
$\nabla \tilde V(x) = \xi(x) + \tau |\xi(x)|^{q-2}A(x)\xi(x)$. We compute
\begin{align*}
 &\quad  \sl^q(u) - |\partial(\phi + \tau \sl^q/q)|^q(u) = \int_{\Tz^d}
  \left(|\nabla V|^q - |\nabla \tilde V|^q \right) {\rm d}
  u \\
  &\quad =\int_{\Tz^d}
  \left(|\xi|^q -  |\xi+ \tau |\xi|^{q-2}A\xi|^q\right) {\rm d}
    u\\
  &\quad  \leq \int_{\Tz^d}
   \left(\left(\big|\xi+ \tau |\xi|^{q-2}A\xi\big|+\big| \tau |\xi |^{q-2}A\xi \big|
    \right)^q -|\xi+ \tau |\xi|^{q-2}A\xi|^q \right) {\rm d}
    u \\
  &\quad =   \int_{\Tz^d} \sum_{k=1}^{\infty}\binom{q}{k}\big|\xi+
    \tau |\xi|^{q-2}A\xi \big|^{q-k} \big| \tau  |\xi|^{q-2}A\xi \big|^k  {\rm d}
    u \\
  &\quad \leq \tau   \sum_{k=1}^{\infty}\binom{q}{k}\|\xi+
    \tau |\xi|^{q-2}A\xi   \|_{L^{\infty}(\Tz^d)}^{q-k} \| |\xi|^{q-2}A\xi \|^k_{L^\infty(\Tz^d)} \\
  &\quad \leq \tau \left(\|\xi+  \tau |\xi|^{q-2}A\xi   \|_{L^{\infty}(\Tz^d)} + \| |\xi|^{q-2}A\xi \|_{L^\infty(\Tz^d)}  \right)^q.
\end{align*}
The one-sided Taylor-expansion condition~\eqref{eq:taylor} hence follows with the choice 
$$g(\tau) =  \tau  \left(\|\nabla V\|_{L^{\infty}(\Tz^d)} +2  \|\nabla V\|^{q-1}_{L^{\infty}(\Tz^d)} \| {\rm
    D}^2V\|_{L^{\infty}(\Tz^d)}    \right)^q.$$
  By applying Theorem \ref{thm:main3} we obtain the following.

\begin{proposition}
  Assume $V\in C^3(\Tz^d)$ and $u^0\in  
\mathcal P(\Tz^d)$. Let  
  $\{0=t_0^n<t_1^n<\dots<t^{n}_{N^n}=T\}$ be a sequence of partitions
  with $\tau^n:=\max (t^n_i-t^n_{i-1}) \to 0$ as $n\to
  \infty$ and $(\tau^n_i - \tau^n_{i-1})^+/\tau^n_{i-1}\leq \haz C
  \tau^n$ for $i=1,\dots,N^n$ . Moreover, let $u^n_i \in M_G(\tau^n_i,u^n_{i-1})$ for $i=1,
  \dots, N^n$ and $\phi$ defined in~\eqref{eq:ff}. Then, up to a not relabeled subsequence, we have that
  $\ove u^n(t) \sto u(t)$, where $u \in AC^p([0,T];\mathcal
  P(\Tz^d))$ satisfies $u(0)=u_0$ and 
the nonlinear transport equation
$$ \partial_t u - {\rm div}\left(  u | \nabla  V|^{q-2} \nabla  V \right)=0\quad \text{in} \ \ \mathcal D'(\Tz^d
\times (0,T)).$$
\end{proposition}

%%%%%%%%%%%%%%%%%%%%%%%%%%%%%%%%%%%%%%%%%%%%%%%%%%%%%%%%%%%%%%%%%%%%%%%%%%%%%%%%%

\section*{Acknowledgements}  
This research is supported by the Austrian Science Fund (FWF) projects
F\,65, W\,1245,  I\,4354, and P\,32788 and by the OeAD-WTZ project CZ 01/2021.

%%%%%%%%%%%%%%%%%%%%%%%%%%%%%%%%%%%%%%%%

\end{document}